\makeatletter
\let\ORIlabel\label
\let\ORIrefstepcounter\refstepcounter
\AddToHook{package/hyperref/before}{%
	\let\label\ORIlabel
	\let\refstepcounter\ORIrefstepcounter
}
\makeatother

\documentclass[onefignum,onetabnum]{siamart220329}

\usepackage{amsmath}
\usepackage{amssymb}
\usepackage{amsfonts}
\usepackage{nccmath}

\usepackage{mathtools}
\usepackage{graphicx}
\graphicspath{ {./images/} }
\usepackage{subcaption}
\usepackage{color}
\usepackage{xcolor}
\usepackage{float}

\usepackage{balance}
\crefformat{equation}{(#1)}
\usepackage{bbm}
\usepackage{bm}
\usepackage{comment}
\usepackage{arcs}

\usepackage{resizegather}

\usepackage{ifthen}

\usepackage{acronym}
\newacro{lqr}[LQR]{Linear Quadratic Regulators}
\newacro{slqr}[SLQR]{Structured Linear Quadratic Regulators}
\newacro{olqr}[OLQR]{Output-feedback Linear Quadratic Regulators}
\newacro{lqg}[LQG]{Linear Quadratic Gaussian}
\newacro{dare}[DARE]{Discrete-time Algebraic Riccati Equation}
\newacro{ouralgo}[RNPO]{Riemannian Newton-type Policy Optimization}
\newacro{PO}[PO]{Policy Optimization}
\newacro{pg}[PG]{Projected Gradient}
\newacro{mdp}[MDP]{Markov Decision Process}

\usepackage{textcomp}
\def\BibTeX{{\rm B\kern-.05em{\sc i\kern-.025em b}\kern-.08em
		T\kern-.1667em\lower.7ex\hbox{E}\kern-.125emX}}

\usepackage{paralist}

\def\x{{X}}
\def\s{{\Lambda}}
\def\r{{\bm r}}
\def\w{{W}}
\def\u{{U}}
\def\H{{\mathcal{H}}}

\def\F{\mathcal{F}}
\def\B{\mathcal{B}}

\def\R{\mathbb{R}}

\def\Amatrices{\mathbb{R}^{n\times n}}
\def\Bmatrices{\mathbb{R}^{n\times m}}

\def\Hmatrices{\mathbb{R}^{n\times d}}
\def\Kmatrices{\mathbb{R}^{m\times n}}

\def\stableK{\mathcal{S}}

\newcommand{\tr}[1]{\ensuremath{\mathrm{tr}\left[ #1 \right]}}

\usepackage{tikz}

\usepackage{amsmath}
\usepackage{amssymb}
\usepackage{nccmath}
\usepackage{mathtools}

\theoremstyle{plain}
\newtheorem{assumption}[theorem]{Assumption}

\usepackage{acronym}
\newacro{lqr}[LQR]{Linear Quadratic Regulators}
\newacro{lqg}[LQG]{Linear Quadratic Gaussian}
\newacro{leqg}[LEQG]{Linear Exponential Quadratic}
\newacro{po}[PO]{Policy Optimization}
\newacro{mdp}[MDP]{Markov Decision Process}
\newacro{cvar}[CVaR]{Conditional Value at Risk}
\newacro{acoe}[ACOE]{Average Cost Optimality Equation}
\newacro{acoi}[ACOI]{Average Cost Optimality Inequality}

\DeclareMathOperator{\E}{\mathrm{\mathbb{E}}}
\def\F{\mathcal{F}}
\renewcommand{\P}{\operatorname{{\mathbb{P}}}}

\usepackage{relsize}
\usepackage{exscale}

\def\B{\mathcal{B}}

\def\R{\mathbb{R}}

\usepackage{lipsum}
\usepackage{amsfonts}
\usepackage{graphicx}
\usepackage{epstopdf}
\usepackage{algorithmic}
\ifpdf
  \DeclareGraphicsExtensions{.eps,.pdf,.png,.jpg}
\else
  \DeclareGraphicsExtensions{.eps}
\fi

\newsiamremark{remark}{Remark}
\newsiamremark{hypothesis}{Hypothesis}
\crefname{hypothesis}{Hypothesis}{Hypotheses}
\newsiamthm{claim}{Claim}

\headers{Ergodic-Risk Criterion for Policy Optimization}{Shahriar Talebi and Na Li}

\title{Ergodic-Risk Criterion for \\ Stochastically Stabilizing Policy Optimization\thanks{A preliminary version of this work, limited to a specific setting and without the main analysis, has been submitted for presentation at the 2025 American Control Conference in Denver, CO.
\funding{NSF AI institute 2112085.}}}

\author{Shahriar Talebi\thanks{UCLA Samueli Mechanical \& Aerospace Engineering Department, University of California, Los Angeles; and School of Engineering and Applied Sciences, Harvard University
  (\email{s.talebi@ucla.edu}, \url{https://shahriarta.github.io}).}
\and Na Li\thanks{School of Engineering and Applied Sciences, Harvard University 
  (\email{nali@seas.harvard.edu}, \url{https://nali.seas.harvard.edu}).}}

\usepackage{amsopn}

\makeatletter
\newcommand*{\addFileDependency}[1]{%
  \typeout{(#1)}%
  \@addtofilelist{#1}%
  \IfFileExists{#1}{}{\typeout{No file #1.}}%
}
\makeatother

\ifpdf
\hypersetup{
  pdftitle={Ergodic-risk Criterion for Stochastically Stabilizing Policy Optimization},
  pdfauthor={Shahriar Talebi and Na Li}
}
\fi

\begin{document}

\maketitle

\begin{abstract}
This paper introduces \textit{ergodic-risk criteria}, which capture long-term cumulative risks associated with controlled Markov chains through probabilistic limit theorems—in contrast to existing methods that require assumptions of either finite hitting time, finite state/action space, or exponentiation necessitating light-tailed distributions.
Using tailored Functional Central Limit Theorems (FCLT), we demonstrate that the time-correlated terms in the ergodic-risk criteria converge under uniform ergodicity and establish conditions for the convergence of these criteria in non-stationary general-state Markov chains involving heavy-tailed distributions.
For quadratic risk functionals on stochastic linear systems, 
in addition to internal stability, this requires the (possibly heavy-tailed) process noise to have \emph{only a finite fourth moment}.
After quantifying cumulative uncertainties in risk functionals that account for extreme deviations, these ergodic-risk criteria are then incorporated into policy optimizations, thereby extending the standard average optimal synthesis to a risk-sensitive framework.
Finally, by establishing the strong duality of the constrained policy optimization, we propose a primal-dual algorithm that optimizes average performance while ensuring that certain risks associated with these ergodic-risk criteria are constrained.
Our risk-sensitive framework offers a theoretically guaranteed policy iteration for the long-term risk-sensitive control of processes involving heavy-tailed noise, which is shown to be effective through several simulations.
\end{abstract}

\begin{keywords}
  Ergodic-risk Criterion; Risk-aware; Linear Quadratic Regulator (LQR); Uniformly Ergodic Chains; Constrained Policy Optimization
\end{keywords}

\begin{AMS}
  93E20
\end{AMS}

\section{Introduction}
\label{sec:intro}

Stochastic control systems yield effective policies for decision-making in uncertain settings like financial markets \cite{rockafellar_optimization_2000}, robotics \cite{majumdar_how_2020}, and healthcare \cite{eichler_risks_2013}. However, optimizing typical average performance often fails to account for risks associated with rare but impactful events, leading to catastrophic consequences. As such, incorporating risk measures become vital in such decision making problems for balancing the performance with resilience to rare events.

While robust control frameworks, such as the mixed $\mathcal{H}_2$-$\mathcal{H}_\infty$ in \cite{zhang_policy_2021}, focus on incorporating the worst-case scenario performance as a constraint, they can be overly conservative when those unlikely events are (possibly) unbounded. Risk-aware approaches \cite{whittle_risk-sensitive_1981,borkar_controlled_1990,borkar_risk-constrained_2014,chow_risk-constrained_2017,sopasakis_risk-averse_2019}, on the other hand, offer a (probabilistic) compromise by building on available stochastic priors to manage both risk and performance simultaneously, offering a better balance by considering both risk and average performance. However, these measures are often deployed on either finite-horizon processes, finite state-action domains, and/or under finite hitting time assumptions \cite{borkar_risk-constrained_2014,chow_risk-constrained_2017,sopasakis_risk-averse_2019}. While avoiding complications regarding limiting probabilities, these assumptions may not fully capture \emph{long-term risk} associated with the stochastic behaviors, especially in \emph{unbounded general-state non-stationary} Markov chains.

Among these, the folklore risk-sensitive framework by Whittle \cite{whittle_risk-sensitive_1981}, aka Linear Exponential Quadratic Gaussian (LEQG), handles the general unbounded nonstationary setting--which (in certain parameter regimes) can be interpreted as optimizing a specific mixture of the average performance and its higher moments (e.g. variance). However, the Gaussian noise is critical for the exponentiation to be well-defined, and thus does not handle cases with \textit{heavy-tailed noise} distributions for capturing rare events; also see the excellent recent survey \cite{biswas_ergodic_2023} expanding on the exponentiation technique\footnote{It is often referred to as \textit{Ergodic risk-sensitive control}; not to be confused with the \textit{ergodic-risk criteria} proposed in this work.}. This motivated \cite{tsiamis_risk-constrained_2020} to introduce a framework for constraining the uncertainty in ``only the state-related portion'' of the finite-horizon \ac{lqr} cost, and \emph{not} the input signal. This is then extended to infinite-horizon through policy optimization techniques \cite{zhao_global_2023}.

In this paper, we introduce the novel \emph{ergodic-risk criteria} to address long-term cumulative risks in stochastic systems by capturing \emph{extreme deviations} from mean performance. This formulation, grounded in functional limit theorems for uniformly ergodic chains, extends risk-constrained control to unbounded, nonstationary Markov chains with (possibly) \emph{heavy-tailed noise}. For quadratic risk functionals, we require the process noise $\w_t$ having bounded moments \textit{only up to the fourth order}, prototyping a framework for risk-aware decision-making in these settings. By incorporating these criteria into constrained optimization frameworks, we formalize the ergodic-risk Constrained Optimal Control Problem (COCP), which balances the average performance with the long-term cumulative risks. 
Our contributions and structure of this paper are organized as follows:
\begin{enumerate}[i)]
    \item {Introducing the Ergodic-risk COCP} (\Cref{sec:probSetup}): We introduce a new COCP using the ergodic-risk criteria, a probabilistic criterion that quantifies long-term cumulative risks, thus accounting for extreme deviations from mean performance. This enables us to address long-term risk in non-stationary processes, previously excluded from the literature (see e.g. \cite{kishida_risk-aware_2023}).

    \item {Convergence Analysis via Uniform Ergodicity} (\Cref{sec:convergence-ergodic-risk}):
    The proposed ergodic-risk criteria are shown to exist for non-stationary processes under uniform ergodicity, ensuring the proposed Ergodic-risk COCP is well-defined.
    Using tailored Functional Limit Theorems, we demonstrate that time-correlated terms in the ergodic-risk criteria converge under uniform ergodicity. For quadratic risk functionals on linear dynamics, uniform ergodicity is ensured under internal stability as long as process noise has a \textit{finite fourth moment}--even if it is heavy-tailed with unbounded higher moments.

    \item {Ergodic-risk Constrained Policy Optimization} (\Cref{sec:quadratic-ergodic-risk}): We explicitly analyze the ergodic-risk criterion for quadratic risk functionals and integrate them into a policy optimization framework. This enables long-term risk-aware control for linear stochastic systems with quadratic risk functional, prototyping risk-sensitive control for unbounded state, non-stationary processes with \emph{heavy-tailed noise distribution}.

    \item {Strong Duality and Primal-Dual Optimization Algorithm} (\Cref{sec:primal-dual}): By establishing a strong duality, we propose a primal-dual method that solves the resulting ergodic-risk COCP by optimizing the system’s average performance and satisfying an ergodic-risk criterion. In particular, this method solves the optimal \ac{lqr} control problem with constraints on the asymptotic conditional variance (as an ergodic-risk), with convergence guarantees and validations through simulations. 
\end{enumerate}
Finally, our ergodic-risk analysis in \Cref{sec:convergence-ergodic-risk} is fundamentally different than that of \cite{tsiamis_risk-constrained_2020,zhao_global_2023} and the proposed ergodic-risk framework (\Cref{thm:C-infty-N-convergence}) generalizes to provably capture the long-term risk associated with heavy-tailed noise for any quadratic functional of \emph{both} the states and the inputs. However similar to \cite{tsiamis_risk-constrained_2020,zhao_global_2023}, the strong duality and the resulting primal-dual algorithm in \Cref{sec:primal-dual} are established only for the case in which the risk functional does not explicitly depend on the inputs.

\section{Problem Setup}
\label{sec:probSetup}

Consider the discrete-time stochastic linear system,
\begin{equation}\label{eq:dynamics}
    \x_{t+1} = A \x_{t} + B \u_{t} + H \w_{t+1}
\end{equation}
where $A \in \Amatrices$, $B \in \Bmatrices$, and $H \in \Hmatrices$ are system parameters as real matrices, $\x_t$ and $\u_t$ denote the state and input vectors, respectively. Also, $\w_t$ denotes the mean-zero process noise and $\x_0$ is the initial state vectors that are sampled independently from probability distributions $\P_\w$ and $\P_0$ with covariances $\Sigma_W$ and $\Sigma_0$, respectively.

\textbf{Notation:} We denote the mean of initial distribution by $m_0 = \E(\x_0)$. We denote the history of state-input trajectory up to time $t$ by $\H_t = \{\x_j,\u_j\}_{j=0}^{t}$ and $\F_t = \sigma(\H_t)$ denotes the $\sigma$-algebra generated by $\H_t$. 
Note that each $\w_t$ is measurable with respect to $\F_{t}$, which we denote by $\w_t \in\F_{t}$.
Let $\F_{-1}$ be the trivial $\sigma$-algebra, then at each time $t\geq0$, we apply an \textit{admissible} input signal $\u_t \in \mathcal{L}_2(\F_{t-1}\vee \sigma\{\x_t\})$ (i.e. square-integrable, measurable with respect to the smallest $\sigma$-algebra containing both $\F_{t-1}$ and $\sigma\{\x_t\}$) and measure the next stochastic state vector $\x_{t+1}$. A \textit{policy} $\pi=\{\pi_t\}_t$ is a collection of measurable mappings that generate an input sequence $\{\u_t = \pi_t(\H_{t-1},\x_t)\}_t$ such $\u_t$ is admissible for all time $t$. A policy is \textit{Markov} if it does not depend on the history $\H_{t-1}$ for every $t$, and called \textit{stationary} if it is independent of time $t$. So, the optimal LQR policy is a stationary Markov policy, and in this work we also restrict ourselves to the class of \emph{affine} stationary Markov policies. We refer to the closed-loop dynamics under a policy $\pi$ as $\mathbf{\Phi}^\pi$ representing the controlled Markov chain.

\subsection{Ergodic-risk Criteria for Markov Chains}
We propose a risk criterion that captures the long-term uncertainty by accumulating step-wise uncertainties as the system evolves. Since states $\{\x_t\}$ of the Markov chain are observed iteratively, it is natural to consider the uncertainty at each stage and accumulate these contributions over time to characterize the overall risk. This leads to a cumulative uncertainty variable, which converges to a limiting value if properly normalized.
 
To formalize this, let us consider any measurable functional of choice 
\[g:\R^n\times\R^m \mapsto \R,\]
called ``risk functional'',
for example a quadratic/affine function in $\x_t$ and/or $\u_t$ which evaluates the performance of each sample path (possibly different than the performance cost). At each time $t-1$, we have access to information in $\F_{t-1}$, so the risk factor at time $t$ is the ``uncertain component'' of the risk functional $g(\x_t,\u_t)$.
This motivates the following definition
\begin{equation}\label{eq:def-C-t}
    C_t \coloneqq 
    g(\x_{t},\u_{t}) -\E[g(\x_{t},\u_{t}) | \F_{t-1}], \quad \text{for $t\geq 0$,}
\end{equation}
 capturing the uncertainty in $g(\x_{t},\u_{t})$ relative to that past information\footnote{In particular, the quadratic choice $g(\x_{t},\u_{t}) = \x_t^\intercal Q\x_t$ captures the state-related uncertainty in the original cost that we aim to minimize; this special term is considered in \cite{tsiamis_risk-constrained_2020} for quantifying ``a risk constraint'' for the uncertainty in the LQR cost.}---see \Cref{fig:one-step-risk}.
To account for the long-term risk behavior, especially in non-stationary systems with heavy-tailed noise, we define the \textit{ergodic-risk criterion} $C_\infty$ as the limit of the normalized cumulative uncertainty:
\begin{equation}\label{eq:def-C-infty}\textstyle
    \frac{1}{\sqrt{t}} S_t \coloneqq \frac{1}{\sqrt{t}} \sum_{s=0}^t C_s \xrightarrow{d\;} C_\infty, \quad \text{as } t\to\infty.
\end{equation}
We also we define the \emph{asymptotic variance} $\gamma_C^2$ as the limit
\begin{equation}\label{eq:asymp-var-def}
     \E[S_t^2/t] \xrightarrow{} \gamma_C^2, \quad \text{as }  t\to \infty.
\end{equation}
Finally, consider the \emph{asymptotic conditional variance} $\gamma_N^2$ defined as the limit
\begin{equation}\label{eq:cond-var-def}\textstyle
     \frac{1}{t} N_t\coloneqq \frac{1}{t} \sum_{s=1}^{t} \E[C_s^2|\F_{s-1}] \xrightarrow{a.s.} \gamma_N^2, \quad \text{as }  t\to \infty.
\end{equation}
The ergodic-risk criterion $C_\infty$, asymptotic variance $\gamma_C^2$, and asymptotic conditional variance $\gamma_N^2$ are tightly related as summarized in the following remark.

\begin{remark}\label{rem:asymp-var-conditional-var}
    The variance of ergodic-risk criterion $C_\infty$ coincides with the asymptotic variance $\gamma_C^2$, whenever they exist 
    (see \Cref{sec:convergence-ergodic-risk} for the existence). Additionally,
    $\gamma_N^2$ serves as an ``estimate'' of $\gamma_C^2$ (whenever exists) in the following sense: By Doob's decomposition we can show that
\(S_t^2 = M_t + N_t\)
where $(M_t,\F_t)$ is a martingale and $(N_t,\F_t)$ is predictable defined as
\begin{equation*}
    N_t - N_{t-1} = \E[S_t^2 - S_{t-1}^2|\F_{t-1}] = \E[(S_t- S_{t-1})^2|\F_{t-1}]
    = \E[C_t^2|\F_{t-1}],
\end{equation*}
where the second equality follows because $(S_t,\F_t)$ is a martingale.
Thus,
\(S_t^2  = M_t + \sum_{s=1}^{t} \E [C_s^2| \F_{s-1}].\)
So, $N_t$ can be interpreted as an intrinsic measure of time for the martingale $S_t$. Also, it can be interpreted as ``the amount of information'' contained in the past history of the process, related to a standard Fisher information \cite[p. 54]{hall_martingale_1980}.
\end{remark}

\subsection{Ergodic-risk Constrained Optimal Control Problem}
Given an input signal $\u = \{\u_t\}_{t=0}^T$, we define the cumulative cost
\begin{align*}\textstyle
    J_{T}(\u) = \sum_{t=0}^T \x_t^\intercal Q \x_t +  \u_t^\intercal R \u_t,
\end{align*}
with $Q \succeq 0$ and
$R\succ 0$ being positive semidefinite and positive definite matrices,
respectively. Conventionally, the infinite-horizon \acf{lqr} problem is to design a sequence of admissible inputs $\u=\{\u_t\}_0^\infty$ that minimizes 
\[\textstyle J(\u) =  \limsup_{T\to\infty} \frac{1}{T} \E [J_T(\u)],\]
subject to dynamics in \cref{eq:dynamics}.

Now, given a coherent risk measure $\varrho:(\mathbb{R},\mathcal{L}) \to \mathbb{R}\cup\{+\infty\}$ from Lebesgue measurable functions to extended reals, we can introduce the ergodic-risk constraint as $\varrho(C_\infty)\leq \varrho_0$ where  $\varrho_0$ is a given risk level--often arising from some initial choice of admissible input sequence.  
However, directly constraining a risk measure on $C_\infty$ leads to a complex, nonlinear constraint on the feedback policy parameters, which in turn hinders efficient policy synthesis.

Nonetheless, as we will establish in \Cref{sec:convergence-ergodic-risk}, $C_\infty$ exists under certain conditions and follows a \emph{normal} distribution. Therefore, a reasonable choice of a coherent risk measure $\varrho$ on $C_\infty$ (such as Conditional Value-at-Risk CVaR$_{\alpha}$ and the Entropic Value-at-Risk EVaR$_{\alpha}$ on $C_\infty$ at a level $\alpha$), essentially reduces to a linear functional of $\gamma_C^2$. Additionally, owing to the relations established in \Cref{rem:asymp-var-conditional-var},
the variance of $C_\infty$ is approximated by the asymptotic conditional variance $\gamma_N^2$. So, $\gamma_N^2$ can be interpreted as an ``estimate'' of $\gamma_C^2$ which has a more tractable expression.

Therefore, for efficient policy synthesis, we instead constrain $\gamma_N^2$ directly, and thus pose the following Ergodic-risk Constrained Optimal Control Problem (COCP):
\begin{align}\label{eq:ergodic-risk-COCP}
   \displaystyle \min \;& J(U) \\
   \text{s.t.} \;& \x_{t+1} = A \x_t + B \u_t + H \w_{t+1}, \quad \forall t \geq0, \nonumber\\
    &\gamma_N^2 \leq \bar\beta, \; \{\u_t\}_{t\geq 0}  \textit{ is an admissible sequence,} \nonumber
\end{align}
for some constant $\bar\beta$ encapsulating the risk level--often chosen from a baseline policy.

As we will see in \Cref{sec:quadratic-ergodic-risk}, this leads to constrained optimization problems where strong duality holds, which enables us to design an efficient primal-dual method for policy optimization with convergence guarantees to the optimal ergodic-risk sensitive policy.

\section{Existence of the Ergodic-risk Criteria}\label{sec:convergence-ergodic-risk}
Before attempting to solve the proposed Ergodic-risk COCP, one need to ensure the well-posedness of this problem. Aside from feasibility of admissible input sequence and the risk constraints, most importantly one need to ensure that these ergodic-risk criteria do exist.

Under minimal assumptions (\Cref{lem:MDS}), it is not hard to show that $\{C_t, \F_t\}_0^\infty$ defined in \cref{eq:def-C-t} is a \emph{Martingale Difference Sequence} (MDS); meaning that for all $t$: i) $C_t$ is $\F_{t}$-adapted (i.e. each $C_t$ is measurable with respect to $\F_{t}$),
ii) $\E |C_t| < \infty$, and iii) $\E[C_t|\F_{t-1}] = 0$ almost surely.
However, establishing the convergences in \cref{eq:def-C-infty,eq:asymp-var-def,eq:cond-var-def} are not possible by directly applying Martingale limit theorems due to the strong correlations in these limiting terms. Instead, we need to obtain much stronger limiting theorems (\Cref{thm:tailored-clt}), building on the uniform ergodicity (\Cref{def:uniform-ergodic}) of the underlying process established in \Cref{thm:V-ergodic}.

Before digression to this general ergodic theory, we discuss the parameterization of stabilizing policies and formalize the claim on the existence of these risk criteria, establishing well-definiteness of the proposed Ergodic-risk COCP for quadratic risk functionals.
Hereafter, we focus on affine ``policies'': 
\[\pi:x\mapsto Kx+\ell\] 
for some matrix parameter $K$ and vector $\ell$, so that each input is designed to be $\u_{k-1} = K \x_{k-1} + \ell$. We define the set of (Schur) stabilizing policies as 
\[\mathcal{S} \coloneqq \left\{K \in \R^{m \times n}\;:\; A_K \coloneqq A+BK \text{ is Schur stable}\right\},\]
i.e., where $A+BK$ has spectral radius less than 1.
We also refer to $\pi = (K,\ell) \in \mathcal{S}\times \R^m$ as the policy without ambiguity. For $\mathcal{S}$ to be non-empty, we consider the following minimal assumption:
\begin{assumption}\label{assmp:stability}
    The pair $(A,B)$ is stabilizable.
\end{assumption}
Also, for each $K\in\mathcal{S}$, let $\Sigma_K$ be the unique positive definite solution to the following Lyapunov equation:
\begin{equation}\label{eq:Sigma-K}
    \Sigma_K = A_K \Sigma_K A_K^\intercal + H\Sigma_W H^\intercal.
\end{equation}

We also require the process noise $\w_t$ and initial states $\x_0$ to be uncorrelated across time, i.e.,
\(\x_0 \perp \w_t, \text{ and } \w_t \perp \w_s, \; \forall t,s \geq 1.\)
In particular, for each $t\geq0$, the process noise $\w_t$ is independent of $\F_{t-1}$, but (by definition) measurable with respect to $\F_t$.
For simplicity, we pose the following stronger assumption.
\begin{assumption}\label{assmp:noise}
    The sequence $\{\w_t\}$ consists of i.i.d. samples of a common zero-mean probability measure $\P_w$ that is non-singular with respect to Lebesgue measure $\lambda$ on $\R^d$, has a non-trivial density, and has finite second moment (with covariance $\Sigma_W \succ0$). %
\end{assumption}

Herein, we focus on quadratic risk functionals:
\begin{equation}\label{eq:def-quadratic-risk-functional}
g(x,u) = x^\intercal Q^c x + u^\intercal R^c u
\end{equation}
for some $Q^c,R^c \succeq 0$. The case of linear risk functional, follows similarly where it requires the noise to have only finite second moment; and so is left to the reader.
The next result establishes the ergodic-risk criteria $C_\infty$, $\gamma_C^2$, and $\gamma_N^2$  in \cref{eq:def-C-infty,eq:asymp-var-def,eq:cond-var-def} are indeed well-defined. 

\begin{theorem}\label{thm:C-infty-N-convergence}
    Suppose Assumptions \ref{assmp:stability} and \ref{assmp:noise}  hold and consider the chain $\mathbf{\Phi}^\pi$ for any policy $(K,\ell) \in \mathcal{S}\times\R^n$ that is stabilizing and $(A_K,H)$ is controllable. Consider the risk functional $g$ in \cref{eq:def-quadratic-risk-functional} and the corresponding $C_t$ in \cref{eq:def-C-t}.
    If the probability measure $\P_w$ generating the noise process $\{W_t\}$ has finite fourth moment, then the followings hold:
    \begin{enumerate}[(i)]
    \item The asymptotic variance $\gamma_C^2$ defined in \cref{eq:asymp-var-def} is well-defined, non-negative, finite, and is given by
    \[\textstyle\gamma_C^2 \coloneqq \sum_{j=1}^d \lambda_j^2 \gamma_{v_j}^2,\]
    where $M = \sum_j \lambda_j v_j v_j^\intercal$ is the eigenvalue decomposition of $M\coloneqq Q_K^c -A_K^\intercal Q_K^c A_K$ with $Q_K^c = Q^c + K^\intercal R^c K$, and $\gamma_{v_j}^2$ is the asymptotic variance given by
    \[\textstyle\gamma_{v_j}^2 = \sum_{k=-\infty}^\infty \Big( \E\left[(v_j^\intercal Y_0 Y_k^\intercal v_j)^2\right] - (v_j^\intercal \Sigma_K v_j)^2 \Big),\]
     where $\{Y_t: t\in\mathbb{Z}\}$ is the stationary process given by $Y_t = \sum_{n = 0}^\infty A_K^n H \w_{t-n}$ with $\{\w_t\}$ as i.i.d. samples of the same probability measure $\P_w$.   
    
    \item As $t\to\infty$ we obtain that
    \[\textstyle\frac{1}{\sqrt{t}} S_t \coloneqq \frac{1}{\sqrt{t}} \sum_{s=1}^t C_s \xrightarrow{d\;} C_\infty \sim \mathcal{N}(0, \gamma_C^2),\]
    if $\gamma_C^2>0$; otherwise $S_t/\sqrt{t}\xrightarrow{a.s.} 0$.
    
    \item As $t\to\infty$ we obtain that $S_t/t \xrightarrow{a.s.} 0,$ and 
    \(\frac{1}{t} \sum_{s=1}^{t} \E [C_s^2| \F_{s-1}] \xrightarrow{a.s.} \gamma_N^2,\)
    where 
    \begin{equation*}
        \hspace{-0.7cm}\gamma_N^2 = 4\tr{Q_K^c H \Sigma_W H^\intercal Q_K^c (\Sigma_K-H \Sigma_W H^\intercal +\bar{x}_K \bar{x}_K^\intercal)} 
        + 4 M_3(Q_K^c)^\intercal Q_K^c \bar{x}_K 
        + m_4(Q_K^c),
    \end{equation*}
    with 
    \begin{equation*}
    \begin{aligned}
        M_3(Q_K^c) &\coloneq \E\big[H\w_{1}\tr{Q_K^c H (\w_{1} \w_{1}^\intercal - \Sigma_W )H^\intercal}\big],\\
        m_4(Q_K^c) &\coloneq \E\left[\tr{Q_K^c H (\w_{1} \w_{1}^\intercal - \Sigma_W )H^\intercal}^2\right].
    \end{aligned}
    \end{equation*}
\end{enumerate}
\end{theorem}
Under the hypothesis of \Cref{thm:C-infty-N-convergence}, by Tower property, we can also deduce that 
\(\frac{1}{t} \sum_{s=1}^{t} \E [C_s^2] \xrightarrow{a.s.} \gamma_N^2.\)

In the rest of this section, we first discuss the necessity of relying on the ergodic theory (\Cref{sec:necessity-ergodic-theory}). Then, provide explicit conditions to ensure uniform egodicity the underlying process (\Cref{thm:V-ergodic} in \Cref{sec:ergodic-LTI}). This enables us to establish a tailored functional limit theorem (\Cref{thm:tailored-clt} in \Cref{sec:tailored-FCLT}), providing us with the theory to finally prove \Cref{thm:C-infty-N-convergence} in \Cref{sec:proofs-of-CLT}. 

\begin{remark}
    Despite the exposition of \Cref{thm:C-infty-N-convergence} for clarity, the implications of forthcoming analysis are far beyond ``quadratic'' risk functionals. For instance, one can directly apply the same result to argue about the ergodic-risk criterion arising from the ``Huber loss'' as risk functional which will be more sensitive to outliers compared to the quadratic one. In particular, the same analysis as in \Cref{thm:V-ergodic} carries through for any $g$ as long as $g^2(x) \leq \|x\|^4 + c$ for some constant $c$.
\end{remark}

\subsection{Setting up the Existence Analysis}\label{sec:necessity-ergodic-theory}
To see why such ergodic theory is required, let us first study properties of $\{C_t, \F_t\}_0^\infty$ with minimal assumptions.
\begin{figure}[!pt]
    \centering
    \begin{tikzpicture}[scale=0.8, transform shape] 
    \draw[thick,->] (-0.5,0) -- (5,0) node[right] {$\mathcal{L}^2(\mathcal{F}_{t-1})$};
    \draw[thick,->] (0,-0.5) -- (0,2.5) node[above right] {\hspace{-1cm}$\mathcal{L}^2(\mathcal{F}_{t-1})^\perp \subset \mathcal{L}^2(\mathcal{F}_{t})$};
    
    \coordinate (O) at (0,0); %
    \coordinate (V) at (3,2); %
    
    \draw[thick,->,line width=1.5pt,>=latex] (O) -- (V) node[above right] {$g(\x_{t},\u_{t})$};
    
    \draw[dashed] (V) -- (3,0) node[below] {\color{blue}$\E[g(\x_{t},\u_{t}) | \F_{t-1}]$};
    \draw[dashed] (V) -- (0,2) node[left] {\color{red}$C_t$};
    
    \draw (3.2,0) -- ++(0,0.2) -- ++(-0.2,0); %
    \draw (0,2.2) -- ++(0.2,0) -- ++(0, -0.2); %
    
    \draw[thick,<-,blue,line width=1.5pt,>=latex] (3,0) -- (0,0);  %
    \draw[thick,<-,red,line width=1.5pt,>=latex] (0,2) -- (0,0);  %
\end{tikzpicture}
    \caption{ The conditional expectation $\E[g(\x_{t},\u_{t}) | \F_{t-1}]$ is the orthogonal projection (in blue) of the risk functional $g(\x_{t},\u_{t})$ onto $\mathcal{L}^2(\mathcal{F}_{t-1})$, i.e. its best estimate by the information up to time $t-1$, i.e. the solution to {\(\arg\min_{\hat g \in \mathcal{L}^2(\mathcal{F}_{t-1})}\|g - \hat g\|_{\mathcal{L}^2(\mathcal{F}_{t-1})} = \sqrt{\E[(g - \hat g)^2]}.\)} 
    The random variable $C_t$ (in red) then retains the ``uncertain part'' of the risk functional at time $t$.}
    \label{fig:one-step-risk}
    \vspace{-.7cm}
\end{figure}

\begin{lemma}\label{lem:MDS}
    Suppose   
    \(|g(x,u)| \leq M(\|x\|^{p} + \|u\|^{p}),\; \forall (x,u) \in\R^n\times\R^m\), and
    for some non-negative constants $M,p \geq 1$.
    Given an affine stationary Markov policy for dynamics \cref{eq:dynamics}, if the process noise $\{\w_t\}$ and initial condition has finite moments up to order $p$, then $\{C_t,\F_{t}\}_0^\infty$ is a Martingale Difference Sequence (MDS); i.e. for all $t$: i) $C_t$ is $\F_{t}$-adapted (i.e. each $C_t$ is measurable with respect to $\F_{t}$),
    ii) $\E |C_t| < \infty$, and iii) $\E[C_t|\F_{t-1}] = 0$ almost surely.
\end{lemma}
Even though, we have shown that $C_t$ is a MDS, the well-definiteness of the limiting quantity in \cref{eq:def-C-infty} requires careful analysis of convergence in distribution. Note that the summands are highly correlated through the dynamics in \cref{eq:dynamics} and thus vanilla Central Limit Theorem does not apply as it requires independent summands. One can instead apply various extended versions of CLT for martingales; e.g. Martingale Central Limit Theorem in \cite[Theorem 5.1]{komorowski_central_2012}. However, even such extension is not useful here because the requirements translate to such strong stability conditions on \cref{eq:dynamics} that is not feasible\footnote{Unless the noise process $\{W_t\}$ is eventually vanishing, which is not the point of interest in this work.}.

Let us define the following processes that become particularly relevant when the risk function is quadratic; namely, the running average and the running covariance:
\begin{align}\label{eq:def-s-t-Gamma-t}
    \s_t \coloneqq \sum_{s=1}^t (\x_s - 
\bar{x}_K),\;  \Gamma_t \coloneqq \sum_{s=1}^t (\x_s - \bar{x}_K) (\x_s -\bar{x}_K)^\intercal, \text{ and } 
    \bar{x}_K \coloneqq \sum_{\tau =0}^{\infty} A_K^\tau B \ell.
\end{align}
Now, we can easily argue about the expectation of running average and covariance quantities defined in \cref{eq:def-s-t-Gamma-t} as follows.

\begin{lemma}\label{lem:ave-conv}
    Under Assumptions \ref{assmp:stability} and \ref{assmp:noise} , for any stabilizing affine stationary Markov policy $\pi \in \mathcal{S}\times \R^m$, we have the following limits as $t\to\infty$:
    \(
        \E[\x_t] \to \bar{x}_K, \; \E[\s_t/t] \to 0,\; \E[\Lambda_t/\sqrt{t}] \to 0, \text{ and } \;
         \E[\Gamma_t/t] \to \Sigma_K,
    \)
    where $\Sigma_K$ is defined in \cref{eq:Sigma-K}.
    Furthermore, we have the following convergence in probability:
    \(\s_t/t \xrightarrow{p} 0, \text{ as } t\to\infty,\)
    and the boundedness in probability: there exists absolute constants $\bar c_1,\bar c_2$ such that for any (large enough) $m$,
    \(\P(\|\x_t\|\geq m) \leq \bar c_1/(m-\bar c_2\|\ell\|)^2, \quad \forall t\geq0.\)
\end{lemma}

This allows us to reason about the first and second moment of the process, and also conclude convergence in probability of terms like
\( \|\x_t\|^2/a_t \xrightarrow{p} 0,\)
for any (scalar deterministic) sequence $a_t\to+\infty$ as $t\to \infty$.
However, it still does not provide much more for the convergence of $C_\infty$ in distribution. In particular, we need much more sophisticated tools to argue about convergence of terms like $\s_t/\sqrt{t}$, $\Gamma_t/\sqrt{t}$, or even $\Gamma_t/t$ as $t\to \infty$ (beyond their expectation).

Next, we build on another type of extensions to CLT known as ``Functional Central Limit Theorem'' that extends the Martingale CLT to Markov chains, connecting to their ``stochastic stability'' properties. For that, we need to establish, in the next section, the so-called ``uniform ergodicity'' of the underlying process as a stochastic stability notion that allows for such Functional CLT to hold. This analysis builds on the ergodic theory \cite{meyn_markov_2009}.

\subsection{$V$-uniform Ergodicity of LTI Systems}\label{sec:ergodic-LTI}
Before proceeding further in this section, we need to introduce the following notation for a general-state \ac{mdp}:

Consider the chain $\mathbf{\Phi}=\{\x_t\}_0^\infty$ that evolves on the Borel space $(\R^n, \B(\R^n))$ according to the linear state space model \cref{eq:dynamics} for some admissible sequence $\{\u_t\}$. Because $\{\w_t\}$ is i.i.d, then this chain has the following stationary transition kernel
\[P_u(x,E) \coloneqq \P \{Ax + B u + H \w_1 \in E\},\]
for any $E \in \B(\R^n)$. We think of a transition probability kernel as
\begin{enumerate}[i)]
    \item $\forall E \in \mathcal{B}(\mathrm{X})$, $P_u(.,E)$ is a non-negative measurable function on $\R^n$, and
    \item $\forall x \in \R^n$, $P_u(x,.)$ is a probability measure on $\mathcal{B}(\R^n)$.
\end{enumerate}
Also, $P_u$ acts on bounded measurable functions $f$ on $\R^n$ as 
\[\textstyle P_u(x,f) \coloneqq P_u f(x) \coloneqq \int_{\R^n} P_u(x, dy) f(y),\] 
and acts on $\sigma$-finite measures $\mu$ on $\mathcal{B}(\R^n)$ as 
\(\textstyle  \mu P_u(E) = \int_{\R^n} \mu(dx) P_u(x,E).\)
Finally, for any stationary Markov policy $\pi$, we denote the MDP resulting from adopting policy $\pi$ by $\mathbf{\Phi}^\pi=\{\x_t\}$ with transition kernel $P_\pi$.
Also, we inductively define the $n$-step transition kernel 
\(\textstyle P_{\pi}^n(x,E) \coloneqq \int_{\R^n} P_{\pi(x)}(x,dy)P_{\pi}^{n-1}(y,E), \; n\geq2,\)
for any $x\in \R^n, E\in\B(\R^n)$, with $P_\pi(x,E) = P_{\pi(x)}(x,E)$. In particular, for every $x\in\R^n$ and $E\in\B(\R^n)$,
\(P_\pi(x,E)=  \P \{Ax + B \pi(x) + \w_1 \in E\} = \P_x \{\x_1 \in E, \} \)
where $\x_1 \in \mathbf{\Phi}^\pi$.

We will establish that, under certain conditions, the closed-loop dynamics \cref{eq:dynamics} under a stabilizing policy has the following properties:
\begin{definition}\label{def:main}
    A chain $\mathbf{\Phi}=\{\x_t\}$ is called
    \begin{enumerate}[(i)]
        \item \emph{Lebesgue-irreducible} if for any measurable set $E \in \B(\R^n)$ with positive Lebesgue measure $\lambda(E)>0$,
        \(\P_x \{\mathbf{\Phi} \text{ ever enters } E\} >0, \; \forall x \in \R^n;\)

        \item \emph{Harris recurrent} if for any $E \in \B(\R^n)$ with $\lambda(E)>0$,
        \(\P_x \{\mathbf{\Phi}\in E ~~i.o.\} =1, \; \forall x \in E,\)
        where $\{\mathbf{\Phi}\in E ~~i.o.\}\coloneqq \bigcap_{k=1}^\infty \bigcup_{t=k}^\infty \{\x_t \in E\}$;

        \item a \emph{T-chain}\footnote{This definition is slightly stronger than the ``T-chain defined in \cite{meyn_policy_1997}'', but it simplifies our exposition here.} if every compact set $E\in \B(\R^n)$ is ``small'': i.e., there exists $m>0$ and some nontrivial measure $\nu_m$ on  $\mathcal{B}(\R^n)$ such that for all $x \in E$ and all $F \in \mathcal{B}(\R^n)$,
\(P^m(x,F) = \P_x \{\x_m \in F \} \geq \nu_m(F);\)
        
        \item \emph{positive} if it admits a unique invariant probability measure $\sigma$, i.e.,
        \(\sigma P(E) = \sigma(E),\quad  \forall E \in \B(\R^n);\)

        \item  \emph{(strongly-)aperiodic} if there exists a $\nu_1$-small set with positive $\nu_1$-measure, i.e., there exists some set $E\in \B(\R^n)$ and some nontrivial measure $\nu_1$ on  $\mathcal{B}(\R^n)$ such that $\nu_1(E) >0$ and for all $x \in E$ and all $F \in \mathcal{B}(\R^n)$,
\(P(x,F) \geq \nu_1(F).\)
    \end{enumerate}
\end{definition}

Despite these properties, we still need to further establish a particular uniform ergodicity for this closed-loop stochastic dynamics, connecting the convergence of the induced transition kernel to the boundedness of \textit{noise moments}. 
To explain further, the standard ergodic theory establishes conditions for existence of the following limit:
\(\lim_{t\to \infty} \E_x[f(\x_t)] = \sigma(f) \coloneqq \int f d\sigma,\)
for every initial condition and every ``bounded function'' $f$ on $\R^n$, where $\sigma$ is an invariant probability measure, i.e. $\sigma(E) = \sigma P_u(E), \;\forall E \in \B(\R^n)$.
To avoid this boundedness assumption (which is too restrictive for our purposes), we instead consider the convergence in the $V$-norm defined as follows:

\begin{definition}\label{def:uniform-ergodic}
Consider two transition kernels $P_1, P_2$ and consider a positive function $V:\R^n\to[1,\infty)$. Define the $V$-norm distance:
\[\textstyle ||| P_1-P_2 |||_V \coloneqq \sup_{x \in \R^n} {\|P_1(x,.) - P_2(x,.)\|_V}\big/{V(x)},\]
where the $V$-norm of any signed measure $\nu$ is defined as:
\[\textstyle \|\nu\|_V \coloneqq \sup_{g:|g|\leq V} \left|\nu(g)\coloneqq \int g d\nu\right|.\]
The chain $\mathbf{\Phi}$ is called \textit{$V$-uniformly ergodic} if 
\(|||P^n - \mathbf{1} \otimes \sigma|||_V \to 0 \; \text{ as } n\to \infty,\)
where $\sigma$ is an invariant probability measure and $\otimes$ denotes the outer product of a function and a measure.    
\end{definition}

Finally, consider the chain $\mathbf{\Phi}^\pi = \{\x_t\}$ representing the closed loop dynamics \cref{eq:dynamics} under a policy $\pi = (K,\ell)$ starting from some $\x_0 = x$.
Now, we establish the properties defined in \Cref{def:main,def:uniform-ergodic} for $\mathbf{\Phi}^\pi$ in the following theorem.

\begin{theorem}\label{thm:V-ergodic}
    Suppose Assumptions \ref{assmp:stability} and \ref{assmp:noise}  holds. Then, for any $\pi \in \mathcal{S}\times \R^m$ such that $(A_K,H)$ is controllable,  Chain $\mathbf{\Phi}^\pi=\{\x_t\}$ is 
    \begin{enumerate}[(i)]
        \item Lebesgue-irreducible, aperiodic T-chain;

        \item positive, Harris recurrent with a unique invariant probability measure $\sigma_\pi$ which is the law of $\x_\infty$ defined as
        \[\textstyle \x_\infty \coloneqq \bar{x}_K + \sum_{t=0}^\infty A_K^t H W_t, \]
        where $\bar{x}_K \coloneqq \sum_{\tau =0}^{\infty} A_K^\tau B \ell$;
        
        \item $V$-uniformly ergodic with $V(x)=\|x - \bar{x}_K\|_M^2+1$ for any $M \succ 0$, where $\|x\|_M^2 \coloneqq x^\intercal M x$, if the noise $\{\w_t\}$ has finite second moment;

        \item $V$-uniformly ergodic with $V(x)=\|x - \bar{x}_K\|_M^4+1$ for any $M \succ 0$, if in addition the noise $\{\w_t\}$ has finite fourth moment. 
    \end{enumerate}
\end{theorem}

As we will see later, for convergence of terms like $\Gamma_t/\sqrt{t}$ we require the chain $\mathbf{\Phi}^\pi$ to be $V$-uniformly ergodic for $V(x) =\|x\|^4+1$. This is the main reason for the requirement that noise process $\{\w_t\}$ must have finite fourth moment. 

\subsection{The tailored Functional LLN and CLT}\label{sec:tailored-FCLT}

Next, building on the uniform ergodicity, we obtain tailored Functional Law of Large Numbers (LLN) and Central Limit Theorem (CLT) with explicit formulas. It enables us to argue about the convergence of linear functionals on $\Gamma_t/t$ and even $\Gamma_t/\sqrt{t}$ with $\Gamma_t$ defined in \cref{eq:def-s-t-Gamma-t}--which was not possible before. Further, it will be used for quantifying ergodic-risks with a quadratic functional $g$.

\begin{theorem}\label{thm:tailored-clt}
        Suppose Assumptions \ref{assmp:stability} and \ref{assmp:noise}  holds and consider the chain $\mathbf{\Phi}^\pi$ for any policy $\pi=(K,\ell) \in \mathcal{S}\times\R^n$ that is stabilizing and $(A_K,H)$ is controllable.
    
    (i) If the noise process $\{W_t\}$ has finite second moment, then 
    \[\textstyle \frac{1}{t} \s_t \xrightarrow{a.s.} 0, \; \text{ and } \frac{1}{\sqrt{t}} \s_t \xrightarrow{d\;} \mathcal{N}(0,\Sigma_\Gamma(0)),\]
    where $\bar{x}_s \coloneqq \sum_{\tau =0}^{s} A_K^\tau B \ell$ and 
     \[\Sigma_\Gamma(\omega) = (I-e^{i\omega} A_K)^{-1} H \Sigma_W H^{\intercal}\left(I-e^{-i\omega} A_K^{\intercal}\right)^{-1}.\]

    (ii) If the noise process $\{W_t\}$ has finite fourth moment, then for any constant matrix $M$
     \[\textstyle \frac{1}{t}\tr{M\Gamma_t} \xrightarrow{a.s.} \tr{M \Sigma_K},\]
     where $\Sigma_K$ solves 
    \(\Sigma_K = A_K \Sigma_K A_K^\intercal + H\Sigma_W H^\intercal.\)
    In this case,
    \begin{equation}\label{eq:gamma-M}\textstyle 
        \frac{1}{\sqrt{t}}\tr{M(\Gamma_t- t \Sigma_K)} \xrightarrow{d\;} \mathcal{N}\left(0,\gamma_M^2 \coloneqq \sum_{j=1}^d \lambda_j^2 \gamma_{v_j}^2\right),
    \end{equation}
    where $\bar M = \sum_j \lambda_j v_j v_j^\intercal$ is the symmetric part of $M$, $\gamma_{v_j}^2$ is the asymptotic variance given by
    \[\textstyle \gamma_{v_j}^2 = \sum_{k=-\infty}^\infty \left( \E\left[(v_j^\intercal Y_0 Y_k^\intercal v_j)^2\right] - (v_j^\intercal \Sigma_K v_j)^2 \right),\]
     and $\{Y_t: t\in\mathbb{Z}\}$ is the stationary process given by $Y_t = \sum_{n = 0}^\infty A_K^n H \w_{t-n}$ with $\{\w_t\}$ as i.i.d. samples of the same probability measure $\P_w$. If, in addition, $\P_w$ is Gaussian then it admits the integral form
    \(\textstyle \gamma_M^2 = \frac{1}{\pi} \int_{-\pi}^\pi \sum_{j=1}^d \lambda_j^2|v_j^\intercal \Sigma_\Gamma(\omega) v_j|^2 d\omega. \)
\end{theorem} 
\noindent Note that in the Gaussian case, $\gamma_M^2$ is directly linked to the $\mathcal{H}_2$-norm of the closed-loop system; which may not generally hold otherwise.

\subsection{Proofs of \Cref{thm:V-ergodic,thm:tailored-clt}, and the resulting \Cref{thm:C-infty-N-convergence}}\label{sec:proofs-of-CLT}

\begin{proof}[Proof of \Cref{thm:V-ergodic}]
    Claim (i) follows from \cite[Proposition 6.3.5]{meyn_markov_2009} whenever $(A_K, H)$ is controllable. To argue claim (ii), by \cite[$f$-norm Ergodic Theorem 14.0.1]{meyn_markov_2009} it suffices to show the so-called ``geometric drift condition'' holds: i.e., for all $ x \in \R^n$,
\[\Delta V (x) \coloneqq P V(x) - V(x) \leq -\beta V(x) +b I_C(x),\]
 for some $V:\mathbb{R}^n \mapsto[1,\infty)$, a constant $\beta>0$, and a small set $C$ (see \Cref{def:main}.(iii)), where $P V(x) \coloneqq \E\left[V(\x_{t+1} )|\x_t = x\right]$.
    Additionally, for Claims (iii) and (iv), by \cite[Theorem 16.0.1]{meyn_markov_2009} it suffices to show that the geometric drift condition holds for each $V$ defined in (iii) and (iv), respectively; Also, by equivalence of norms (in finite dimensional spaces), it further suffices to do so only for the specific matrix $M$ solving the following Lyapunov equation 
    \[M = A_K^\intercal M A_K + Q,\]
    for some $Q \succ I$. 
    Note that, as $A_K$ is Schur stable, there always exists a unique positive definite $M \succ Q$ solving this Lyapunov equation.
    Also, noting that \Cref{assmp:noise} requires the noise $\{\w_t\}$ has finite second moment, the quadratic case of $V(x)=\|x\|_M^2+1$ is simpler and follows similarly as Claim (iv). 
    Finally, note that by \Cref{lem:ave-conv}, $\E[\x_t]\to \bar{x}_K$ for any $\ell \in \R^m$. Therefore, it suffices to only show Claim (iv) for the centered process $\{\x_t -\bar{x}_K\}$.
    
    Therefore, we focus on showing Claim (iv) for $V(x) = \|x\|_M^4 +1$ with $\ell =0$ implying $\bar{x}_K = 0$. 
    We compute that
    \begin{align*}
        V(\x_{t+1}) =& \|A_K \x_t + \w_{t+1}\|_M^4 + 1\\
        =& \left(\|A_K \x_t\|_M^2 + \|\w_{t+1}\|_M^2 + 2 \langle A_K \x_t, \w_{t+1}\rangle \right)^2  +  1 \\
        =& \|A_K \x_t\|_M^4 + \|\w_{t+1}\|_M^4 + 4|\langle A_K \x_t, \w_{t+1}\rangle|^2\\
        &+ 2\|A_K \x_t\|_M^2 \|\w_{t+1}\|_M^2 + 4 \|A_K \x_t\|_M^2 \langle A_K \x_t, \w_{t+1}\rangle\\
        &+ 4 \|\w_{t+1}\|_M^2 \langle A_K \x_t, \w_{t+1}\rangle + 1\\
        \leq & \|A_K \x_t\|_M^4 + \|\w_{t+1}\|_M^4 + 4\| A_K \x_t\|_M^2 \|\w_{t+1}\|_M^2\\
        &+ 2\|A_K \x_t\|_M^2 \|\w_{t+1}\|_M^2 + 4 \|A_K \x_t\|_M^2 \langle A_K \x_t, \w_{t+1}\rangle\\
        &+ 4 \|\w_{t+1}\|_M^3 \|A_K \x_t\|_M + 1,
    \end{align*}
where the last inequality follows by Cauchy-Schwarz and the fact that $M \succ Q\succ I$, and so $\|x\|\leq \|x\|_Q \leq \|x\|_M \leq \lambda_{\max}(M) \|x\|_Q$ for all $x\in\R^n$.
Thus, as $\w_{t+1}$ is independent of $\F_t$, the conditional expectation below simplifies so for any $x\in\R^n$
    \begin{align*}
        \Delta V(x) =& \E[V(\x_{t+1}|\x_t = x)] - V(x)\\
        \leq & \|A_K x\|_M^4 - \|x\|_M^4 + m_4 + 4 m_2 \| A_K x\|_M^2 + 2 m_2 \|A_K x\|_M^2 + 4 m_3 \|A_K x\|_M 
    \end{align*}
where $m_4 = \E \|\w\|_M^4$, $m_3 = \E\|\w\|_M^3$, and $m_2 = \E \|\w\|_M^2$. The Lyapunov equation implies that for any $x\in\R^n$,
\[\|A_K x\|_M^2 - \|x\|_M^2 = -\|x\|_Q^2.\]
Thus,
    \begin{align*}
        \Delta V(x) =& \E[V(\x_{t+1}|\x_t = x)] - V(x)\\
        \leq & -\|x\|_Q^2 (\|x\|_M^2 + \|A_K x\|_M^2 ) +m_4 +2m_3 + (6m_2 + 2m_3)\|A_K x\|_M^2\\
        = & - \|x\|_Q^2 \|x\|_M^2 + (6m_2 + 2m_3 -\|x\|_Q^2)\|A_K x\|_M^2  + m_4 +2m_3 \\
        \leq & - \frac{1}{2\lambda_{\max}(M)^4}V(x)  + (6m_2 + 2m_3 -\|x\|_Q^2)\|A_K x\|_M^2  + m_4 +2m_3 \\
        &+ \frac{1}{2\lambda_{\max}(M)^4} - \frac{1}{2}\|x\|_Q^4
    \end{align*}
    where the first inequality follows by $2\|A_K x\|_M \leq \|A_K x\|_M^2 + 1$, and the last inequality follows by $\|x\|_Q \leq \|x\|_M \leq \lambda_{\max}(M) \|x\|_Q$. Now, choose the set $C$ to be
    \[C \coloneqq \left\{x \;\big|\; \|x\|_Q \leq \max\{ 6m_2 + 2m_3,  1 + 2m_4 +4m_3\} \right\},\]
    and notice that
    \(\Delta V(x) \leq -\frac{1}{2 \lambda_{\max}(M)^4} V(x) + b I_C(x),\)
    where
    \begin{multline*}
        b \coloneqq m_4 +2m_3 + \frac{1}{2}
        + (6m_2 + 2m_3) \lambda_{\max}(M)^2 \max\{ 6m_2 + 2m_3,  1 + 2m_4 +4m_3\}^2, 
    \end{multline*}
    because $\|A_K x\|_M\leq \|x\|_M \leq \lambda_{\max}(M) \|x\|_Q$ and $\|x\|_Q \leq \max\{ 6m_2 + 2m_3,  1 + 2m_4 +4m_3\}$ on $C$.
    Finally, as long as $\pi$ is stabilizing and the noise has finite fourth moment, $b$ is a uniformly bounded constant and $C$ is a compact measurable set because $x \mapsto \|x\|_Q$ is continuous. By the first claim, $\mathbf{\Phi}^\pi$ is a Lebesgue-irreducible strongly aperiodic T-chain, and thus $C$ is small set and thus also a petite set. Therefore, the geometric drift condition holds, and the claim follows by \cite[Theorem 16.0.1]{meyn_markov_2009}. 
    
    Finally, in order to obtain a computable expression for $\sigma_\pi$ in the most general case, one can use the construction of ``minimal sub-invariant measures'' \cite[Theorem 10.4.9]{meyn_markov_2009}. For the particular LTI case, we provide the following argument which is very similar to \cite[Section 10.5.4]{meyn_markov_2009}.
    By \Cref{assmp:noise}, \cref{eq:sys-traj} implies that 
    \[\x_{t+1} \sim A_K^{t+1} \x_0 + \bar{x}_t + \sum_{s =0}^{t} A_K^{s}H\w_s,\]
    i.e. has the same law. But, by \Cref{lem:ave-conv}, we observe that as $A_K$ is stable
    \begin{align*}
        \E[\|\x_{t+1}\|] &\leq \|A_K^{t+1}\|m_0+ (\|\bar{x}_K\| +\E[\|W_1\|]) \sum_{s=0}^t \|A_K^s H\|\\
        &\leq C_K m_0 +  (\|\bar{x}_K\| +\E[\|W_1\|]) \sum_{s=0}^\infty \|A_K^s H\|<\infty,
    \end{align*}
    for some constant $C_K>0$ where the second inequality follows by \cite[Lemma 6]{talebi_data-driven_2023} and the last one is due to finite moment condition on $W_1$. Thus, Fubini-Tonelli's Theorem implies that $\x_\infty$ is absolutely convergent and well-defined with 
    \[\E[\|\x_\infty\|]\leq (\|\ell\| + \E[\|W_1\|]) \sum_{s=0}^\infty \|A_K^s H\|<\infty.\] 
    
    Next we show that $\sigma_\pi$, as the law of $\x_\infty$, is an invariant measure. Consider any bounded continuous function $g$ and any $x_0\in\R^n$, and note that by the Dominated Convergence Theorem
    \(P^t g(x_0) = \E_{x_0}g(\x_t) \to \E g(\x_\infty),\)
    as $t\to \infty$.
    Now, consider the function $Pg$ and note that it is continuous as $g$ is and $P$ represents a T-chain, a specially Feller continuous kernel. So,
    \begin{equation*}
        \sigma_\pi(Pg) = \E[Pg(\x_\infty)] = \lim_{t\to\infty}\E_{x_0} Pg(\x_t) =\lim_{t\to\infty}\E_{x_0} g(\x_{t+1}) = \E[g(\x_\infty)] = \sigma_\pi(g).
    \end{equation*}
    Since $\pi_\sigma$ is determined by its value over continuous bounded functions, it must be the unique invariant probability measure.
\end{proof}

\begin{proof}[Proof of \Cref{thm:tailored-clt}]
    Similar to the argument in the proof of \Cref{lem:ave-conv} for the expectations \cref{eq:exp-Lambda-t} and \cref{eq:expec-Gamma-t}, recall that by exponential stability, both
    \(\sum_{s=1}^t (\bar{x}_{s-1} - \bar{x}_K)\)
    and $\sum_{s=1}A_K^s m_0$ converge exponentially fast as $t\to \infty$.
    Thus, it suffices to show the claim for the centered process $\{\x_t -\bar{x}_K\}$ and thus consider only the policies with $\ell=0$. 
      
    By \Cref{thm:V-ergodic}, for any policy $K\in\mathcal{S}$, the model $\mathbf{\Phi}^\pi = \{\x_t\}$ with $\pi=(K,0)$ is $V$-uniformly ergodic with $V(x) =\|x\|^2 +1$ whenever noise process $\{W_t\}$ has finite second moment, and with $V(x) =\|x\|^4 +1$ whenever it has finite fourth moment. 
    Also recall that,
    \(\sigma_\pi(g) = \int g d \sigma_\pi = \E[g(\x_\infty)].\)
    For any function $g\colon \R^n \to \R$ and invariant probability measure $\sigma_\pi$ under policy $\pi$, define 
    \(\bar{g}(x) \coloneqq g(x) - \sigma_\pi(g)\) and consider the series 
    \(S_t^\pi(\bar{g}) \coloneqq \sum_{k=1}^{t} \bar{g}(\x_k),\)
    where $\mathbf{\Phi}^\pi = \{\x_k\}$ is the Markov chain under policy $\pi$.  
    So, for stabilizing policies, $\mathbf{\Phi}^\pi$ is also a positive Harris chain with invariant probability $\sigma_\pi$. Therefore, by \cite[Theorem 17.0.1]{meyn_markov_2009}, the Functional LLN and CLT hold for any such functional $g$. We now utilize this result to show Claims (i) and (ii).
    
    Claim (i): For the first part, this implies that for any vector $v\neq0 \in \R^n$, the LLN and CLT holds for process $\boldsymbol{Y} = \{v^\intercal \x_t\}$ with $g(x) =x$ because $g^2 \leq V$; i.e. the LLN holds:
    \(
        \frac{1}{t}\sum_{k=1}^{t} v^\intercal \x_k \xrightarrow{a.s.} \sigma_\pi^v(g).
    \)
    But by \Cref{thm:V-ergodic}, $\sigma_\pi^v$ is the law of $v^\intercal \x_\infty$ and thus
    \(\sigma_\pi^v(g) = \E[g(v^\intercal \x_\infty)] = v^\intercal\E[ \x_\infty ] = 0,\)
    by the absolute convergence of series defining $\x_\infty$.
    Also, the CLT holds:    
    \begin{equation}\label{eq:clt-linear}
        \textstyle v^\intercal \s_t/\sqrt{t} = \frac{1}{\sqrt{t}}\sum_{s=1}^{t} v^\intercal \x_s \xrightarrow{d\;} \mathcal{N}(0,\gamma_v^2) 
    \end{equation}
    with the asymptotic variance $\gamma_v^2>0$. Next, by \cite[Theorem 17.6.2]{meyn_markov_2009}, the asymptotic variance $\gamma_v$ is given by
    \begin{equation}\label{eq:gamma-v}
         \gamma_v^2=v^{\intercal}(I-A_K)^{-1} H \Sigma_W H^{\intercal}\left(I-A_K^{\intercal}\right)^{-1} v.
    \end{equation}
    To see this, let $Y_t=v^{\intercal} \x_t$ for $t \in \mathbb{Z}_{+}$. Since $A_K$ is Schur stable, then a unique invariant probability $\sigma_\pi$ exists by \Cref{thm:V-ergodic}, and hence a stationary version of the process $Y_t$ also exists, defined for $t \in \mathbb{Z}$. The stationary process can be realized as
    \(
    Y_t=\sum_{n=0}^{\infty} h_{n} \w_{t-n},
    \)
    where $h_{n}=v^{\intercal} A_K^{n} H$ and $\left\{\w_t\right\}_{t \in \mathbb{Z}}$ are i.i.d. with mean zero and covariance $\Sigma_W=$ $\mathrm{E}\left[\w_1 \w_1^{\intercal}\right]$, which is assumed to be finite. Now if we define the autocovariance sequence 
    \(R(s) = \E_{\sigma_\pi^v}[Y_s Y_0], \quad s\in\mathbb{Z},\)
    then the asymptotic variance is given by
    \(\gamma_v^2 = \sum_{s=-\infty}^\infty R(s) = D(\omega)|_{\omega = 0},\)
    where $D(\omega)$ is the spectral density associated with stationary process $\mathbf{Y}$ and the equality follows by the Fourier series form of $D(\omega)$. But it is easy to show that \cite[pp. 66]{kumar_stochastic_2015}
    \(D(\omega) = F(e^{i\omega})\Sigma_W F(e^{i\omega})^*,\)
    where
    \( F(e^{i\omega}) \coloneqq \sum_{n=0}^\infty h_n e^{in\omega} = v^\intercal (I-e^{i\omega} A_K)^{-1} H.\)
    Thus, evaluating $D(\omega)|_{\omega=0}$  justifies the expression in \cref{eq:gamma-v}. As choice of $v\in\R^n$ is arbitrary, the first claim then follows by Cram\'er-Wold device \cite{durrett_probability_2019}.
   
    Claim (ii): Similarly for the second part, this implies that for any vector $v\neq0 \in \R^n$, the LLN and CLT holds for process $\boldsymbol{Y} = \{v^\intercal \x_t\}$ with $g(x) =\|x\|^2$ because $g^2 \leq V$; i.e. the LLN holds:
    \(
        \frac{1}{t}\sum_{s=1}^{t} v^\intercal \x_s \x_s^\intercal v \xrightarrow{a.s.} \sigma_\pi^v(g).
    \)
    In this case,
    \(\sigma_\pi^v(g) = \E[g(v^\intercal \x_\infty)] = v^\intercal\E[ \x_\infty \x_\infty^\intercal ]v = v^\intercal \Sigma_K v, \)
    where the last inequality follows by the noise assumption \Cref{assmp:noise} and the absolutely convergence series form of Discrete Lyapunov Equation in \cref{eq:Sigma-K}.
    Additionally, the CLT holds:    
    \begin{equation}\label{eq:clt-vxxTvT}\textstyle
        \frac{1}{\sqrt{t}}\sum_{s=1}^{t} v^\intercal (\x_s \x_s^\intercal - \Sigma_K) v \xrightarrow{d\;} \mathcal{N}(0,\gamma_v^2) 
    \end{equation}
    with the asymptotic variance $\gamma_v^2>0$ which in this case,
    considering the stationary process $\mathbf{Y}$, is given by
    \begin{align*}
        \gamma_v^2 = \sum_{s=-\infty}^\infty \E_{\sigma_\pi^v}\left[(Y_0^2 - \sigma_\pi^v(g))(Y_s^2 - \sigma_\pi^v(g))\right]
        = \sum_{s=-\infty}^\infty \E_{\sigma_\pi^v}\left[Y_0^2 Y_s^2\right] - \E_{\sigma_\pi^v}[Y_0^2] \E_{\sigma_\pi^v}[Y_s^2],
    \end{align*}
    where the second equality is because $\mathbf{Y}$ is the stationary process, so \( \sigma_\pi^v(g) = \E_{\sigma_\pi^v}[Y_s^2] = v^\intercal \Sigma_K v\) 
    for any $s \in \mathbb{Z}$. Now, if the noise process $\{\w_t\}$ is Gaussian, then $\mathbf{Y}$ is a Gaussian process and by Isserlis' Theorem we obtain that
    \[\E_{\sigma_\pi^v}\left[Y_0^2 Y_s^2\right] - \E_{\sigma_\pi^v}[Y_0^2] \E_{\sigma_\pi^v}[Y_s^2] = 2 (\E_{\sigma_\pi^v}\left[Y_0 Y_s\right])^2,\]
    and thus
    \(\gamma_v^2 = 2 \sum_{s=-\infty}^\infty R(s)^2 = \frac{1}{\pi} \int_{-\pi}^\pi |D(\omega)|^2 d\omega, \)
    where the last equality follows by Parseval's theorem with $|\cdot|$ denoting the modulus. 
    Finally, for any matrix $M$ and any positive semi-definite matrix $\Gamma$, note that
    \[\tr{M\Gamma} = \tr{\Gamma^{1/2}M\Gamma^{1/2}}= \tr{\Gamma^{1/2}\bar M\Gamma^{1/2}} = \tr{\bar M\Gamma},\]
    where $\bar M$ is the symmetric part of $M$. Therefore, it suffices to show the claim for symmetric $M$ matrices. Consider the eigenvalue decomposition of $M=\sum_{i=1}^d \lambda_i v_i v_i^\intercal$ and notice that  
    \(\Gamma \mapsto \tr{M \Gamma} = \sum_{i=1}^d \lambda_i v_i \Gamma v_i^\intercal\)
    is a continuous mapping.
    By \cref{eq:clt-vxxTvT}, for each $i\in\{0,1,\cdots,d\}$ we have
    \( \frac{1}{\sqrt{t}} v_i (\Gamma_t -t\Sigma_K) v_i^\intercal \xrightarrow{d\;} \bar\Gamma_i\)
    where $\bar\Gamma_i \sim \mathcal{N}(0,\gamma_{v_i}^2)$ with $\gamma_{v_i}^2$ defined in \cref{eq:gamma-v} for each $v_i$. But $v_i \perp v_j$ whenever $i\neq j$ because $M$ is symmetric. This implies that $\{\bar\Gamma_i\}_{i=1}^d$ is a set of mutually uncorrelated Gaussian random variables and thus independent. Therefore, the last claim follows by Continuous Mapping theorem and because
    \(\sum_{i=1}^d \lambda_i \bar\Gamma_i \sim \mathcal{N}\left(0,\sum_{i=1}^d \lambda_i^2 \gamma_{v_i} \right).\)
This completes the proof.
\end{proof}

\begin{proof}[Proof of \Cref{thm:C-infty-N-convergence}] \label{sec:proof-C-infty-N-convergence}
    Consider the process $\{\x_t\}$ with any $K\in\mathcal{S}$ and starting from a fixed $\x_0 = x_0$.
    For simplicity, we define $G_{t} := g(\x_{t},\u_{t}), ~t\geq 0$ and note that $G_0 = x_0^\intercal Q_K^c x_0,$ and
    \(
     G_{t+1} = (A_K\x_t+B\ell+H\w_{t+1})^\intercal Q_K^c (A_K\x_t+B\ell+H\w_{t+1}),
    \)
    for $t\geq 0$. 
     This, together with \Cref{assmp:noise} imply  that
    \begin{multline*}
    \E[G_{t+1}|\F_{t}] = \x_t^\intercal A_K^\intercal Q_K^c A_K\x_t + \tr{Q_K^c H \Sigma_W H^\intercal} 
        + 2\ell^\intercal B^\intercal Q_K^c A_K \x_t + \ell^\intercal B^\intercal Q_K^c B \ell
    \end{multline*}
    where $\Sigma_W$ denotes the covariance of $\w\sim\P_\w$. Thus,
    \begin{multline}\label{eq:C-t-proof}
        C_{t+1} = \x_{t+1}^\intercal Q_K^c \x_{t+1} -\x_t^\intercal A_K^\intercal Q_K^c A_K\x_t \\
        - \tr{Q_K^c H \Sigma_W H^\intercal} - 2\ell^\intercal B^\intercal Q_K^c A_K \x_t - \ell^\intercal B^\intercal Q_K^c B \ell,
    \end{multline}
    for $t \geq 0$. By definitions in \cref{eq:def-s-t-Gamma-t}, we obtain that
    \begin{equation}
        \label{eq:def-r-t-Lmabda-t}
        \begin{aligned}
            \r_t &\coloneqq \sum_{s=1}^t \x_s = \s_t + t \bar{x}_K,\; \text{ and }\;
        Z_t &\coloneqq \sum_{s=1}^t \x_s \x_s^\intercal =\Gamma_t + 2 \bar{x}_K \s_t^\intercal + t \bar{x}_K \bar{x}_K^\intercal.
        \end{aligned}
    \end{equation}
    Recall $S_{t} \coloneqq \sum_{s=1}^{t} C_s$ and therefore, the cyclic property of trace and the definitions $\r_t$ and $Z_t$ imply that
    \begin{multline*}
        S_{t} = \tr{Q_K^c Z_{t}} -\tr{A_K^\intercal Q_K^cA_K(Z_{t-1} + x_0 x_0^\intercal)} 
         - t\, \tr{Q_K^c H \Sigma_W H^\intercal} \\
         - 2\ell^\intercal B^\intercal Q_K^c A_K (\r_{t-1} + x_0)
         - t\, \ell^\intercal B^\intercal Q_K^c B \ell,
    \end{multline*}
    which, using \cref{eq:def-r-t-Lmabda-t}, can be rewritten as
    \begin{align}\label{eq:S-t-expression}
        S_{t} =& \tr{Q_K^c \Gamma_{t}} + 2 \bar{x}_K ^\intercal Q_K^c \s_t + t \bar{x}_K^\intercal Q_K^c \bar{x}_K
        -\tr{A_K^\intercal Q_K^cA_K(\Gamma_{t-1} + x_0 x_0^\intercal)}\\
        &- 2 \bar{x}_K^\intercal A_K^\intercal Q_K^cA_K \s_{t-1}
        - t \bar{x}_K^\intercal A_K^\intercal Q_K^cA_K \bar{x}_K  - t\, \tr{Q_K^c H \Sigma_W H^\intercal}\\
        &- 2\ell^\intercal B^\intercal Q_K^c A_K (\s_{t-1} + (t-1)\bar{x}_K + x_0)
         - t\, \ell^\intercal B^\intercal Q_K^c B \ell,
    \end{align}
    Also, by Lyapunov equation \cref{eq:Sigma-K}, we have
    \begin{equation}\label{eq:trace-lyap-identity}
        \tr{Q_K^c\Sigma_K -A_K^\intercal Q_K^c A_K \Sigma_K - Q_K^c H \Sigma_W H^\intercal} = 0.
    \end{equation}
    Therefore, by applying the LLN in \Cref{thm:tailored-clt} to \cref{eq:S-t-expression} we conclude that as $t\to\infty$, 
    \begin{multline*}
        S_t/t \xrightarrow{a.s.}  
        \bar{x}_K^\intercal Q_K^c \bar{x}_K - \bar{x}_K^\intercal A_K^\intercal Q_K^c A_K \bar{x}_K
        - 2\ell^\intercal B^\intercal Q_K^c A_K \bar{x}_K
        -\ell^\intercal B^\intercal Q_K^c B \ell\\
        = \bar{x}_K^\intercal Q_K^c \bar{x}_K - (A_K \bar{x}_K + B\ell)^\intercal Q_K^c (A_K \bar{x}_K + B\ell)
        =0,
    \end{multline*}
    where the last equality follows by the fact that the stationary average $\bar{x}_K$ must satisfy
    \begin{equation}\label{eq:ell-bar-identity}
        \bar{x}_K = A_K \bar{x}_K + B\ell.
    \end{equation}
    Next, by \cref{eq:trace-lyap-identity} we can rewrite \cref{eq:S-t-expression} also as 
    \begin{multline}\label{eq:S-t-expression-clt}
        S_{t} = \tr{(Q_K^c -A_K^\intercal Q_K^cA_K) (\Gamma_{t}-t\Sigma_K)}  \\
        +\tr{A_K^\intercal Q_K^cA_K((\x_t-\bar{x}_K) (\x_t-\bar{x}_K)^\intercal - x_0 x_0^\intercal)} 
        + 2 \bar{x}_K ^\intercal Q_K^c \s_t - 2 \bar{x}_K^\intercal A_K^\intercal Q_K^cA_K \s_{t-1} \\
        - 2\ell^\intercal B^\intercal Q_K^c A_K \s_{t-1}
        + t \bar{x}_K^\intercal Q_K^c \bar{x}_K - t \bar{x}_K^\intercal A_K^\intercal Q_K^cA_K \bar{x}_K - 2 t \ell^\intercal B^\intercal Q_K^c A_K \bar{x}_K \\
        - t\, \ell^\intercal B^\intercal Q_K^c B \ell  
        + 2\ell^\intercal B^\intercal Q_K^c A_K (\bar{x}_K - x_0)
         \\
         = \tr{(Q_K^c -A_K^\intercal Q_K^cA_K) (\Gamma_{t}-t\Sigma_K)} 
        +\tr{A_K^\intercal Q_K^cA_K((\x_t-\bar{x}_K) (\x_t-\bar{x}_K)^\intercal - x_0 x_0^\intercal)} \\
        + 2 \bar{x}_K ^\intercal Q_K^c (\x_t-\bar{x}_K) 
        + 2\ell^\intercal B^\intercal Q_K^c A_K (\bar{x}_K - x_0)
    \end{multline}
    where the last equality follows by \cref{eq:ell-bar-identity}.
    Next, recall that the CLT in \Cref{thm:tailored-clt} implies the convergence in distribution of $\frac{1}{\sqrt{t}}(\tr{M(\Gamma_{t} - t\Sigma_K)}$ for any constant matrix $M$. Furthermore, \Cref{lem:ave-conv} implies that both $\mathrm{tr}[A_K^\intercal Q_K^c A_K (\x_t-\bar{x}_K)(\x_t-\bar{x}_K)^\intercal]/\sqrt{t}$ and $\bar{x}_K ^\intercal Q_K^c (\x_t-\bar{x}_K)/\sqrt{t}$ converges to zero in probability as $t\to\infty$.
    Now, consider the linear (and thus continuous) mapping $\Gamma \mapsto \tr{(Q_K^c - A_K^\intercal Q_K^c A_K)\Gamma}$ and therefore, by applying Continuous Mapping Theorem to \cref{eq:S-t-expression-clt}, we obtain that $\frac{1}{\sqrt{t}}S_{t}$ converges in distribution to $\mathcal{N}(0,\gamma_C^2)$ whenever $\gamma_C^2>0$, and otherwise converges to zero almost surely.
    Finally, we define the conditional covariance
    \[N_{t}^2 \coloneqq \sum_{s=1}^{t} \E [C_{s+1}^2| \F_{s}],\]
    and show the convergence of $N_t^2/t$.
    For that, \cref{eq:C-t-proof} simplifies to
    \begin{align*}
        C_{t+1} =& 2 ( A_K \x_t + B \ell)^\intercal Q_K^c H \w_{t+1}+ \tr{Q_K^c H (\w_{t+1} \w_{t+1}^\intercal - \Sigma_W )H^\intercal}\\
        =& 2 (A_K (\x_t-\bar{x}_K) + \bar{x}_K)^\intercal Q_K^c H \w_{t+1} + \tr{Q_K^c H (\w_{t+1} \w_{t+1}^\intercal - \Sigma_W )H^\intercal}
    \end{align*}
    where the last equality follows by \cref{eq:ell-bar-identity}.
    Thus
    \begin{align*}
        \E[C_{t+1}^2|\F_{t}] =& 4 (A_K (\x_t-\bar{x}_K) + \bar{x}_K)^\intercal Q_K^c H \Sigma_W H^\intercal Q_K^c (A_K (\x_t-\bar{x}_K) + \bar{x}_K) \\
        &+ 4 (A_K (\x_t-\bar{x}_K) + \bar{x}_K)^\intercal Q_K^c \\
        &\cdot\E[H\w_{t+1}\tr{Q_K^c H (\w_{t+1} \w_{t+1}^\intercal - \Sigma_W )H^\intercal} |\F_{t}] \\
        &+\E[\tr{Q_K^c H (\w_{t+1} \w_{t+1}^\intercal - \Sigma_W )H^\intercal}^2|\F_{t}]\\
        =& 4 (A_K (\x_t-\bar{x}_K) + \bar{x}_K)^\intercal Q_K^c H \Sigma_W H^\intercal Q_K^c (A_K (\x_t-\bar{x}_K) + \bar{x}_K) \\
        &+ 4 (A_K (\x_t-\bar{x}_K) + \bar{x}_K)^\intercal Q_K^c M_3(Q_K^c) + m_4(Q_K^c)\\
        =& 4 \tr{ A_K^\intercal Q_K^c H \Sigma_W H^\intercal Q_K^c A_K (\x_t-\bar{x}_K)(\x_t-\bar{x}_K)^\intercal}  \\
        &+ \bar{x}_K^\intercal Q_K^c H \Sigma_W H^\intercal Q_K^c (8 A_K (\x_t-\bar{x}_K)+ 4\bar{x}_K) \\
        &+ 4 M_3(Q_K^c)^\intercal Q_K^c (A_K (\x_t-\bar{x}_K) + \bar{x}_K)  + m_4(Q_K^c)
    \end{align*}
    where $M_3(Q_K^c)$ and $m_4(Q_K^c)$ are defined in \cref{thm:C-infty-N-convergence} and we dropped the conditionals because $\w_{t+1}$ is independent of $\F_t$. Both quantities in \cref{thm:C-infty-N-convergence} are well-defined by the moment assumption on the noise.  
    Therefore,
    \begin{multline*}
        N_{t}^2 =  4\tr{A_K^\intercal Q_K^c H \Sigma_W H^\intercal Q_K^c A_K \Gamma_t} 
        + \bar{x}_K^\intercal Q_K^c H \Sigma_W H^\intercal Q_K^c (8 A_K \s_t + 4 t \bar{x}_K)\\
        + 4 M_3(Q_K^c)^\intercal Q_K^c (A_K \s_t + t\bar{x}_K) + t\, m_4(Q_K^c).
    \end{multline*}
    Finally, by \Cref{thm:tailored-clt}, we obtain the almost sure convergence: 
    \begin{multline*}
        \frac{1}{t}N_t^2 \xrightarrow{a.s.} 4\tr{A_K^\intercal Q_K^c H \Sigma_W H^\intercal Q_K^c A_K \Sigma_K} 
        + 4 \bar{x}_K^\intercal Q_K^c H \Sigma_W H^\intercal Q_K^c \bar{x}_K \\
        + 4 M_3(Q_K^c)^\intercal Q_K^c \bar{x}_K + m_4(Q_K^c).
    \end{multline*}
    But, by cyclic permutation property of the trace and the Lyapunov equation:
    \begin{equation*}
        \tr{A_K^\intercal Q_K^c H \Sigma_W H^\intercal Q_K^c A_K \Sigma_K} 
        = \tr{Q_K^c H \Sigma_W H^\intercal Q_K^c (\Sigma_K-H \Sigma_W H^\intercal)}.
    \end{equation*}
    Combining the last two equations completes the proof.
\end{proof}

\section{Policy Optimization for Quadratic Ergodic-risk COCP} \label{sec:quadratic-ergodic-risk}
In this section, we show that how the functional limit theorems developed in the previous section can be utilized to quantify limiting risk criteria that are useful in control and decision-making. In particular, such limiting risk criteria can be incorporated as a constraint in the optimal control and reinforcement learning frameworks. 

\subsection{Constrained Policy optimization}\label{sec:CPO}
The policy optimization (PO) approach to control design pivots on a parameterization of the feasible policies for the synthesis problem. 
One can view the \ac{lqr} cost naturally as a map $J(K,\ell) : (K,\ell) \mapsto J(\u=K \x+\ell)$. 
And, for any stabilizing policy $\pi = (K,\ell) \in \mathcal{S}\times \R^m$, we can compute the cost as
\begin{multline*}
    J_T(\pi) = 
    \x_0^\intercal Q_K \x_0 + \sum_{t=1}^T \big[(\x_t -\bar{x}_K)^\intercal Q_K (\x_t -\bar{x}_K) 
    + 2\bar{x}_K^\intercal Q_K (\x_t - \bar{x}_K) + \bar{x}_K^\intercal Q_K \bar{x}_K \big]\\
    =\tr{Q_K\left(\x_0 \x_0^\intercal + \Gamma_T+ 2\s_T\bar{x}_K^\intercal + T\bar{x}_K \bar{x}_K^\intercal\right)}
\end{multline*}
with $\s_t$, $\Gamma_t$, and $\bar{x}_K$ defined in \cref{eq:def-s-t-Gamma-t}, $A_K = A+BK$, $Q_K = Q+K^\intercal R K$, and where we used the cyclic permutation property of trace to obtain the last equality. Now, by \Cref{lem:ave-conv} and \Cref{thm:C-infty-N-convergence} we obtain
\[J(\pi) = \lim_{T\to \infty} E[J_T(\pi)/T] = \tr{Q_K(\Sigma_K + \bar{x}_K \bar{x}_K^\intercal)},\]
with $\Sigma_K$ in \cref{eq:Sigma-K} and $\bar{x}_K$ in \cref{eq:def-s-t-Gamma-t}.

    In the absence of any risk constraint, it is well known that the optimal solution to this problem reduces to solving the \ac{dare}
\begin{equation*}
    P_{\mathrm{LQR}} = A^\intercal P_{\mathrm{LQR}} A + Q
    - A^\intercal P_{\mathrm{LQR}} B (R + B^\intercal P_{\mathrm{LQR}} B)^{-1} B^\intercal P_{\mathrm{LQR}} A, 
\end{equation*}
for the unknown matrix $P_{\mathrm{LQR}}$ that quadratically parameterizes the so-called cost-to-go.
Subsequently, one sets $\u_t^*= K_{\mathrm{LQR}} \x_t$, where the optimal \ac{lqr} gain (policy) $K_{\mathrm{LQR}} \in \Kmatrices$ is given by $K_{\mathrm{LQR}} = - (R + B^\intercal P_{\mathrm{LQR}} B)^{-1} B^\intercal P_{\mathrm{LQR}}  A,$
and the optimal cost $J(\u^*) = \tr{P_{\mathrm{LQR}} H \Sigma_W H^\intercal}$ with $\Sigma_W$ denoting the covariance of $\P_\w$. Note that $\Sigma_0$ will not affect the optimal cost.

Finally, by applying \Cref{thm:C-infty-N-convergence}, the optimization problem in \cref{eq:ergodic-risk-COCP} is well-defined and reduces to:
\begin{align}\label{eq:optimization-reform}
    \displaystyle \min_{(K,\ell) \in \stableK \times \R^m} &\; J(K,\ell)  = \tr{Q_K(\Sigma_K + \bar{x}_K \bar{x}_K^\intercal)} \\
   \text{s.t.}~~  \Sigma_K & = A_K \Sigma_K A_K^\intercal + H \Sigma_W H^\intercal, \nonumber\\
     \gamma_N^2 &\leq \bar\beta \nonumber
\end{align}
with $\gamma_N^2$ characterized in \cref{thm:C-infty-N-convergence}, where 
\begin{equation}\label{eq:bar-ell-inverse-ident}
    \bar{x}_K = (I-A_K)^{-1} B \ell,
\end{equation}
which follows by the identity in \cref{eq:ell-bar-identity}.
So, let us consider the Lagrangian $L:\mathcal{S} \times \R^m \times \R_{\geq 0} \mapsto \R$ defined as
\begin{align*}
    L(K, &\ell, \lambda) \coloneqq \tr{Q_K(\Sigma_K+\bar{x}_K \bar{x}_K^\intercal)} + \lambda \left(\gamma_N^2-\bar\beta  \right) \nonumber\\
     =& \tr{(Q_K + 4\lambda Q_K^c H \Sigma_W H^\intercal Q_K^c)(\Sigma_K+\bar{x}_K \bar{x}_K^\intercal)} + 4 \lambda M_3(Q_K^c)^\intercal Q_K^c \bar{x}_K - \lambda \beta(Q_K^c), 
\end{align*}
where 
\(
    \beta(Q_K^c) \coloneqq 4\tr{(Q_K^c H \Sigma_W H^\intercal)^2 } - m_4(Q_K^c) +\bar\beta,
\)
with $m_4(Q_K^c)$ defined in \cref{thm:C-infty-N-convergence}.

\subsection{Quadratic Ergodic-risk Criteria with $R^c=0$ and $\ell = 0$} \label{subseq:R-c-0}
Hereafter, we assume that the risk measure $g$ does not depend on control input explicitly; i.e. $R^c=0$, and so $Q_K^c = Q^c$ is constant.
In addition to system-theoretic assumptions in \Cref{assmp:stability} that is necessary for feasibility of the optimization (non-empty domain $\stableK$), in the following result we also need controllability of the pair $(A^\intercal, Q^{\frac{1}{2}})$ to guarantee regularity of the Lagrangian, i.e. coercivity of $L(\cdot,\lambda)$ for each $\lambda$.
For simplicity, we consider slightly stronger conditions and obtain the following results similar to \cite[Lemma 1 and 2]{talebi_data-driven_2023}.

\begin{assumption}\label{asmp:Q-H-full-row-rank}
    Assume $Q \succ 0 $ and $H$ is full row rank.
\end{assumption}

\begin{proposition}\label{prop:Lagrangian}
Suppose Assumptions \ref{assmp:stability}, \ref{assmp:noise}, and \ref{asmp:Q-H-full-row-rank} hold.
For each $\lambda\geq0$, consider the Lagrangian $L_\lambda(\cdot) = L(\cdot,\lambda):\stableK \to \R$.
The following statements are true:
\begin{enumerate}[(i)]
    \item $L_\lambda(\cdot)$ and $\gamma_N^2(\cdot)$ are smooth with 
    \(
    \nabla L_\lambda(K) = 2(RK + B^\intercal P_{(K,\lambda)} A_K)\Sigma_K 
    \)
    where $P_{(K,\lambda)}$ is the unique solution to
    \[P_{(K,\lambda)} = A_K^\intercal P_{(K,\lambda)} A_K + Q_K + 4\lambda Q^c H \Sigma_W H^\intercal Q^c.\]
    \item $L_\lambda(\cdot)$ is coercive with compact sublevel sets $\stableK_\alpha$ for each $\alpha>0$.
    \item $L_\lambda(\cdot)$ admits a unique global minimizer $K^*(\lambda) = \arg \min_{K\in\stableK} L_\lambda(K)$ that is stabilizing, given by
    \[K^*(\lambda) = -(R + B^\intercal P_{(K^*(\lambda),\lambda)}B)^{-1}B^\intercal P_{(K^*(\lambda),\lambda)} A,\]
    and $L_\lambda(K^*(\lambda)) = \tr{P_{(K^*(\lambda),\lambda)} H \Sigma_W H^\intercal} + \lambda \beta[Q^c]$.
    \item The restriction $L_\lambda(\cdot) |_{\stableK_\alpha}$ for any (non-empty) sublevel set $\stableK_\alpha$ has Lipschitz continuous gradient, is gradient dominated, and has a quadratic lower model; in particular, for all $K,K'\in\stableK_\alpha$ the following inequalities are true: 
    \(\|\nabla L_\lambda(K) - \nabla L_\lambda(K')\|_F \leq \ell \|K - K'\|_F\)
    and
    \[c_2 \|K-K^*\|_F^2\leq c_1[L_\lambda(K) - L_\lambda(K^*)] \leq \|\nabla L_\lambda(K)\|_F^2,\]
    for some positive constants $\ell(\alpha),c_1(\alpha)$, and $c_2(\alpha)$ that only depend on $\alpha$ and are independent of $K$.
\end{enumerate}

\end{proposition}
\begin{proof}
    The smoothness of $L_\lambda(\cdot)$ and $\gamma_N^2(\cdot)$ follows by smoothness of $\Sigma_K$ in $K$, which together with the expression for $\nabla_K L(K,\lambda)$ are obtained similar to derivations of \cite[Lemma 2]{talebi_data-driven_2023}, with $P_{(K,\lambda)}$ being the unique positive definite solution to the displayed Lyapunov equation whenever $K\in\mathcal{S}$, i.e., it is stabilizing. Also, $Q_K + 4\lambda Q^c H \Sigma_W H^\intercal Q^c \succeq Q$ for each $\lambda \geq 0$, and thus
    \(L(K,\lambda) + \lambda \beta(Q^c) \geq \tr{Q_K \Sigma_K} = \tr{P_{(K,0)} H \Sigma_W H^\intercal},\)
    where the equality follows by the cyclic permutation property of trace with $P_{(K,0)}$ satisfying 
    \(P_{(K,0)} = A_K^\intercal P_{(K,0)} A_K + Q_K.\)
    Now, by the cyclic property of trace, we can show that
    \(L(K,\lambda) + \lambda \beta[Q^c] \geq \tr{P_n H_n \Sigma_W H_n^\intercal},\)
    where $P_n$ is the unique solution to 
    \(P_n = (A_K^n)^\intercal P_n A_K^n + Q_K\)
    which is positive definite as $(A_K,Q^{\frac{1}{2}})$ is observable, and $H_n = [H \;\; A_K H \;\; \cdots \;\; A_K^{n-1} H]$ which is full-rank as $(A_K, H)$ is controllable; thus $H_n \Sigma_W H_n^\intercal$ is positive definite as $\Sigma_W \succ 0$, and $\underline{\lambda}(H_n \Sigma_W H_n^\intercal)\geq \underline{\lambda}(H \Sigma_W H^\intercal)>0$ as $H$ is full row rank. Also, by a similar computation as \cite[Lemma 1]{talebi_data-driven_2023}, we can show that 
    \(\tr{P_n H_n \Sigma_W H_n^\intercal} \geq \underline{\lambda}(H_n \Sigma_W H_n^\intercal) \underline{\lambda}(R)\|K\|_F^2.\)
    Therefore, if $\|K\| \to \infty$, $L(K,\lambda)$ approaches $+\infty$; also, if $K\to \partial \stableK$, $P_n$ approaches $+\infty$, resulting in $L(K,\lambda)$ approaching $+\infty$ because $\underline{\lambda}(H_n \Sigma_W H_n^\intercal)\geq \underline{\lambda}(H \Sigma_W H^\intercal)>0$. So, we conclude the coercivity of $L(\cdot,\lambda)$ and thus, the resulting compact sublevel sets. The rest of the claims follows directly by observing that $P_{(K,\lambda)}$ essentially satisfies the dual Lyapunov equation as that $X_{(L)}$ in \cite[Lemma 2]{talebi_data-driven_2023} where $Q$ is replaced with $Q + 4 \lambda Q^c H \Sigma_W H^\intercal Q^c$. This completes the proof.
\end{proof}

\subsection{Strong Duality}
Now that $K^*(\lambda)$ is uniquely well-defined for each $\lambda \geq 0$, the \emph{dual problem} of \cref{eq:optimization-reform} (with $R^c=0$ and $\ell = 0$) can be written in following forms
\begin{align}\label{eq:min-max-problem}
    \sup_{\lambda\geq0} \min_{K \in \stableK} L(K,\lambda) 
&= \sup_{\lambda\geq0} \tr{P_{(K^*(\lambda),\lambda)} H \Sigma_W H^\intercal} + \lambda \beta[Q^c]\\
&=  \sup_{\lambda\geq0} \tr{Q_{K^*(\lambda)} \Sigma_{K^*(\lambda)}} + \lambda \left(\gamma_N^2(K^*(\lambda))-\bar\beta  \right), \nonumber
\end{align}
with $P_{(K^*(\lambda),\lambda)}$ defined in \Cref{prop:Lagrangian}.
It is standard to assume the primal problem is strictly feasible: 
\begin{assumption}[Slater's Condition]\label{asmp:slater}
    Assume $\bar\beta$ is large enough such that there exists $K\in\stableK$ with $\gamma_N^2(K) < \bar\beta$.
\end{assumption}

\begin{proposition}\label{prop:strong-duality}
        Under Assumptions \ref{assmp:stability}, \ref{assmp:noise}, \ref{asmp:Q-H-full-row-rank}, and \ref{asmp:slater}, strong duality holds for the sup-min problem in \cref{eq:min-max-problem} with a unique solution $(K^*(\lambda^*), \lambda^*)$; i.e. both primal and dual optimization problems are feasible with identical values:
    \[\max_{\lambda\geq0} \min_{K \in \stableK} L(K,\lambda)  = L(K^*(\lambda^*), \lambda^*) =\min_{K \in \stableK} \max_{\lambda\geq0} L(K,\lambda). \]
\end{proposition}
\begin{proof}
Because we know the cost function is globally lower bounded, $\tr{Q_{K} \Sigma_{K}} \geq 0$, then Slater's condition implies feasibility of the dual problem, and that there exists a finite $\lambda_0\geq 0$ such that $\gamma_N^2(K^*(\lambda_0)) \leq \bar\beta$. Now, if we let 
\begin{equation}\label{eq:lambda-star}
\lambda^* = \min\{\lambda_0 \;:\; \lambda_0 \geq 0 \text{ and } \gamma_N^2(K^*(\lambda_0)) \leq \bar\beta\}, 
\end{equation}
then we claim that the pair $(K^*(\lambda^*), \lambda^*)$ is the saddle point of the Lagrangian $L(K,\lambda)$, and therefore obtain the strong duality. For this claim to hold, by the Karush–Kuhn–Tucker (KKT) conditions and \Cref{prop:Lagrangian} (iii), it suffice to show that the complementary slackness holds:
\begin{equation}\label{eq:comp-slack}
    \mathrm{CS}(\lambda^*) \coloneqq \lambda^*(\gamma_N^2(K^*(\lambda^*)) - \bar\beta) = 0.
\end{equation}
But, we know that $P_{(K,\lambda)}$ is real analytic in $(K,\lambda)$ and is positive definite for each $\lambda\geq0$. Thus, $K^*(\lambda)$ as defined in \Cref{prop:Lagrangian} is smooth in $\lambda$. So, $\gamma_N^2 \circ K^*(\cdot)$ is continuous by composition and lower bounded by zero. Therefore, complementary slackness follows because the minimum of a strictly monotone function on a compact set in \cref{eq:lambda-star} is unique and is attained at the boundary. This completes the proof.
\end{proof}

In this particular case where the risk functional does not directly depend on the input and the controller is linear (i.e., $R^c=0$ and $\ell = 0$), by \Cref{prop:Lagrangian} the optimal ergodic-risk policy $K^*(\lambda^*)$ can be interpreted as an optimal LQR controller for the same dynamics but where $Q$ is replaces with $Q+4\lambda^* Q^c H \Sigma_W H^\intercal Q^c$. This implies that the optimal ergodic-risk policy is indeed more conservative in penalizing the state deviations from zero. Nonetheless, it depends on the value of $\lambda^*$ that satisfies the CS condition which must be solved for.
Luckily, the established strong duality enables us to provide an efficient algorithm next. 

\section{The Algorithm} \label{sec:primal-dual}
By establishing a strong duality, we can solve the dual problem without loss of optimality. In particular, using the properties of $L_\lambda(\cdot)$ obtained in \Cref{prop:Lagrangian}, we devise simple primal-dual updates to solve \cref{eq:optimization-reform} by accessing the gradient and constrained violation values, as proposed in \Cref{algo}.
\begin{algorithm}[ht]
\caption{Primal-Dual Ergodic-risk Constrained LQR}
    \begin{algorithmic}[1]
\STATE Set $K_0 \in \stableK_\alpha$ for some $\alpha>0$, $\lambda_0 =1$, and tolerance $\epsilon>0$ and stepsize $\eta_m =(\gamma_N^2(K_0) - \bar\beta)^{-1} (m+1)^{-1/2}$
\FOR{$m=0,1,\cdots, T=\mathcal{O}(\ln(\ln(\epsilon))/\epsilon^2)$}
\WHILE{$\|\nabla_K L_{\lambda_m}(K)\|_F \geq \sqrt{\epsilon}$}
    \STATE $G \gets -(R + B^\intercal P_{K,{\lambda_m}} B)^{-1} \nabla L_{\lambda_m}(K) \Sigma_K^{-1}$
    \STATE $K \gets K + \frac{1}{2} G$ 
\ENDWHILE
\STATE $\lambda_{m+1} \gets \max\left[0, \lambda_{m} + \eta_m (\gamma_N^2(K) - \bar\beta) \right]$
\ENDFOR
\RETURN $(K,\lambda \coloneq \frac{1}{T}\sum_{m=1}^{T} \lambda_m)$
\end{algorithmic}
\label{algo}
\end{algorithm}

We next provide a convergence guarantee by combining recent LQR policy optimization \cite{fazel_global_2018, talebi_policy_2023} with saddle point optimization techniques \cite{nedic_subgradient_2009}, and illustrate the performance of \Cref{algo} through simulations.

\subsection{Convergence Guarantee}\label{subsec:convergence-simulations}
The convergence of \Cref{algo} essentially follows from the results obtain in \Cref{prop:Lagrangian}, combined with the results in \cite{talebi_policy_2023} on convergence of quasi-Newton gradient descent for the LQR problem.

\begin{proposition}
    Under Assumptions \ref{assmp:stability}, \ref{assmp:noise}, \ref{asmp:Q-H-full-row-rank}, and \ref{asmp:slater}, \Cref{algo} converges to an $\epsilon$-solution of max-min optimization \cref{eq:min-max-problem} in $\mathcal{O}(\ln(\ln(\epsilon))/\epsilon^2)$ steps.
\end{proposition}
\begin{proof}
 \Cref{prop:Lagrangian} and the initialization in \Cref{algo} ensure that the premise of \cite[Theorem 4.3]{talebi_policy_2023} is satisfied. Also, $G$ is the Riemannian quasi-Newton direction and the update on $K$ is known as the Hewer's algorithm which is proved to converge at a quadratic rate, as discussed in detail in \cite[Remark 5]{talebi_policy_2023}. And therefore, the inner loop terminates very fast in $\ln(\ln(\epsilon))$ steps and essentially returns an $\epsilon$ accurate estimate of $K^*(\lambda_m)$.
Furthermore, the outer loop is expected to take $\mathcal{O}(1/\epsilon^2)$ steps to obtain $\epsilon$ error on the functional $L(K^*(\lambda), \lambda)$ following standard primal-dual guarantees \cite{nedic_subgradient_2009}. Thus, assuming that Assumptions \ref{assmp:noise}, \ref{assmp:stability}, \ref{asmp:Q-H-full-row-rank}, and \ref{asmp:slater} hold, \Cref{algo} obtains an $\epsilon$ accurate solution to problem \cref{eq:optimization-reform} in $\mathcal{O}(\ln(\ln(\epsilon))/\epsilon^2)$ steps. 
\end{proof}

Instead of employing quasi-Newton updates for the inner-loop in \Cref{algo}, one could opt for pure gradient descent updates; however, this approach would lead to a slower convergence rate.
We could similarly notice that the premises of Theorem 1 in \cite{talebi_data-driven_2023} are also satisfied using properties in \Cref{prop:Lagrangian}. Then, as a special case of \cite[Theorem 1]{talebi_data-driven_2023} where we have accurate gradient information $\nabla L_\lambda(K)$, and constants $\gamma, s$ and $s_0$ approach zero, we obtain the linear convergence of gradient descent to the unique global minimizer $K^*(\lambda)$. The inner loop then only utilizes the gradient information $\nabla_K L_\lambda(K)$, however, it takes $\mathcal{O}(\ln(\epsilon))$ steps to complete. Finally, combining either of these algorithms with gradient estimates obtained from input-output data \cite{fazel_global_2018,talebi_policy_2024} seamlessly provides data-driven algorithms for ergodic risk control with end-to-end sample complexity, that is left to the reader.

\subsection{Simulations}\label{sec:simulation}
Herein, we present two simulations: the first evaluates the performance of the ergodic-risk optimal policy against the LQR optimal policy under heavy-tailed noise, while the second illustrates the convergence behavior of \Cref{algo} on randomly generated problem instances.

\textbf{Simulation 1:} We first show how the ergodic-risk optimal policy behaves compared to the LQR optimal policy. For that, we consider the Grumman X-29 aircraft dynamics which, for the purpose of agile maneuverability, was designed with a high degree of longitudinal static instability. The longitudinal and lateral-directional dynamics matrices in the Normal Digital Up-and-Away (ND-UA) mode are reported in \cite[Tables 13 and 14]{bosworth_linearized_1992}, with fixed discretization step-size 0.05 (also see \cite{talebi_regularizability_2022} for online stabilization of the same plant). The dynamics consists of 4 longitudinal states, 4 lateral-directional states, and 5 inputs, and as it is normalized we use consider $Q = Q_c = I_8$ and $R = R_c = I_5$. 

Next, in order to simulate gust disturbances particularly affecting the unstable longitudinal dynamics, we inject a disturbance of magnitude 20 at every 200 time steps in addition to random variables sampled from Student's t-distribution with parameter $\nu =5$. Note that this noise distribution has finite fourth moment and unbounded fifth moment. Note that the standard risk-sensitive control syntheses (based on exponentiation of cost functional) such as LEQG do not apply as the higher moments of such noise distribution do not exist.

\begin{figure}[pt]
  \centering    \includegraphics[width=0.8\textwidth]{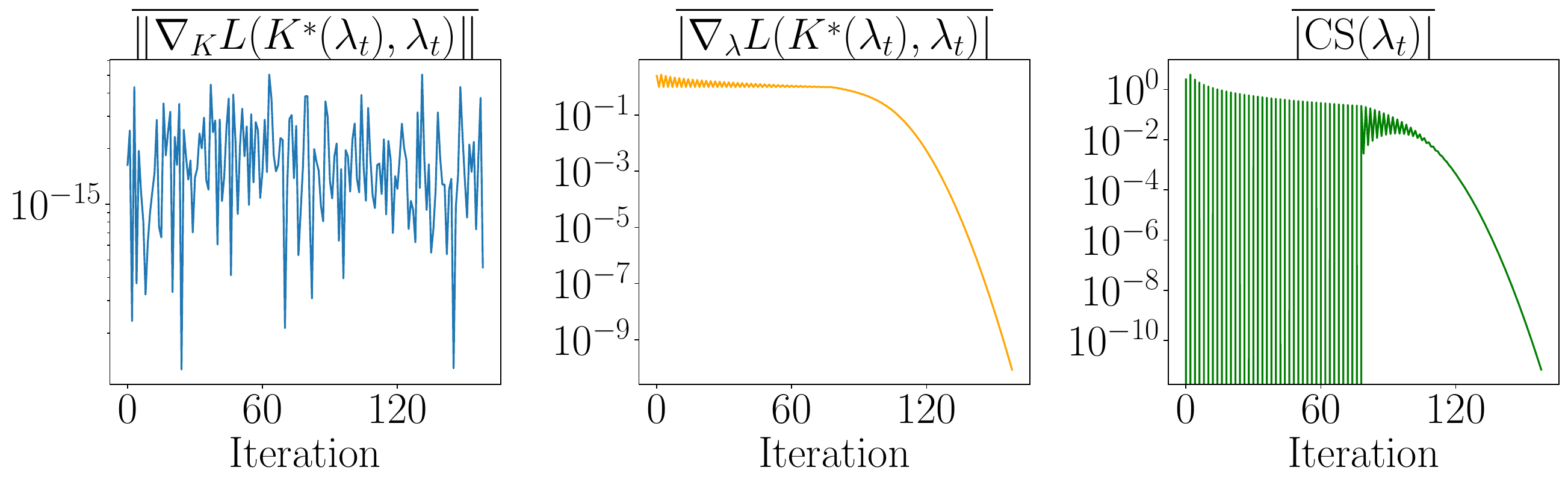}
  \caption{ Convergence of \Cref{algo} to the optimal ergodic-risk policy for Grumman X-29 aircraft dynamics.}
  \label{fig:single-convergence}
  \vspace{-0.7cm}
\end{figure}

Now, the optimal ergodic-risk policy for a risk level of $\bar\beta = 0.8 \gamma_N^2(K_{LQR})$ is obtained via \Cref{algo}, where the errors in KKT conditions are illustrated in \Cref{fig:single-convergence} as performance measures and $\overline{(\cdot)}$ denotes normalization of a sequence with respect to its first element. 
We obtain $J(K_{\mathrm{LQR}}) = 116527$ and $J(K^*) = 116854$. While the ergodic-risk policy $K^*$ is $20\%$ more risk-sensitive (in terms of the asymptotic conditional variance $\gamma_N^2$), it cost only $0.2\%$ more compared with the optimal LQR policy.

\begin{figure}[pt]
  \centering    \includegraphics[width=0.8\textwidth]{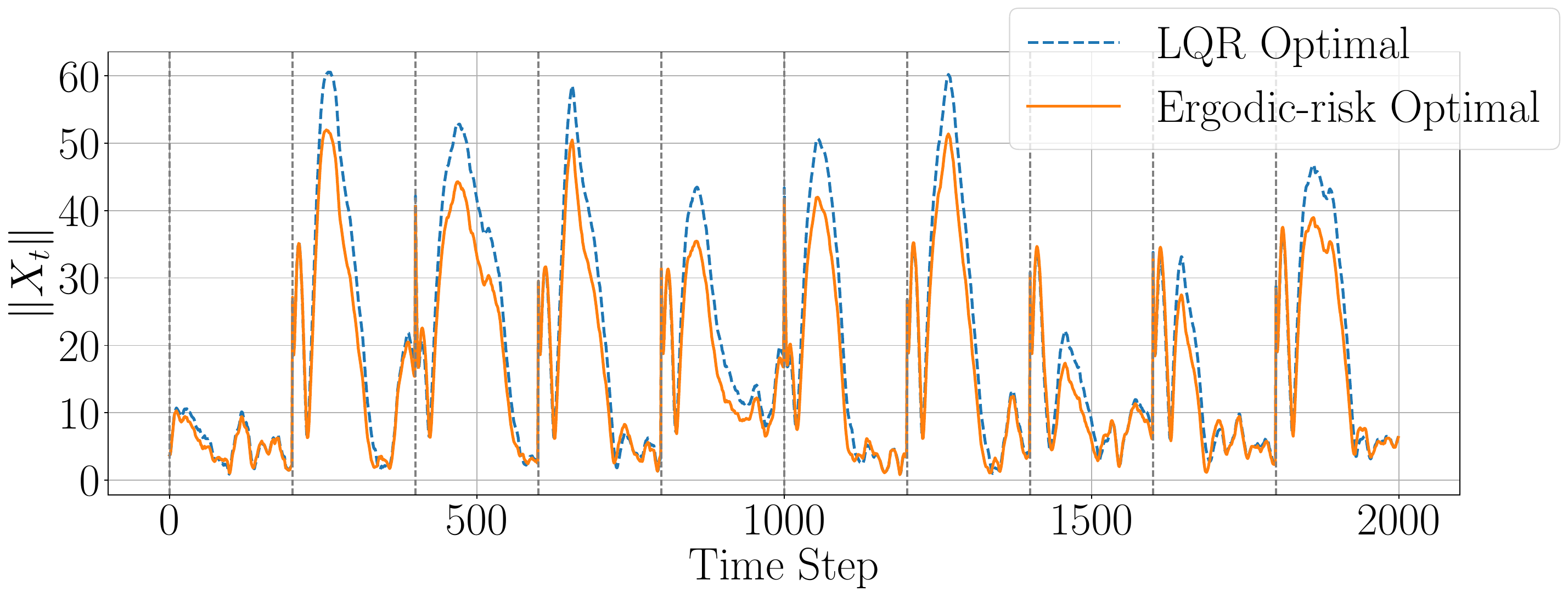}
  \caption{ The optimal ergodic-risk versus the optimal LQR policies for Grumman X-29 aircraft dynamics under simulated gust disturbances at every 200 time steps.}
  \label{fig:ergodic-vs-lqr-rollout}
  \vspace{-0.5cm}
\end{figure}

Next, performance of the optimal ergodic-risk versus the optimal LQR polices are illustrated for a single roll out in \Cref{fig:ergodic-vs-lqr-rollout}.
Evidently, the ergodic-risk policy is more resilient to gust disturbances compared with the LQR policy.
Similar behavior of resilience to large disturbances for such ergodic-risk sensitive policies are also observed in \cite{tsiamis_risk-constrained_2020,zhao_global_2023}.

\begin{figure}[pt]
    \centering
    \includegraphics[width=0.8\textwidth]{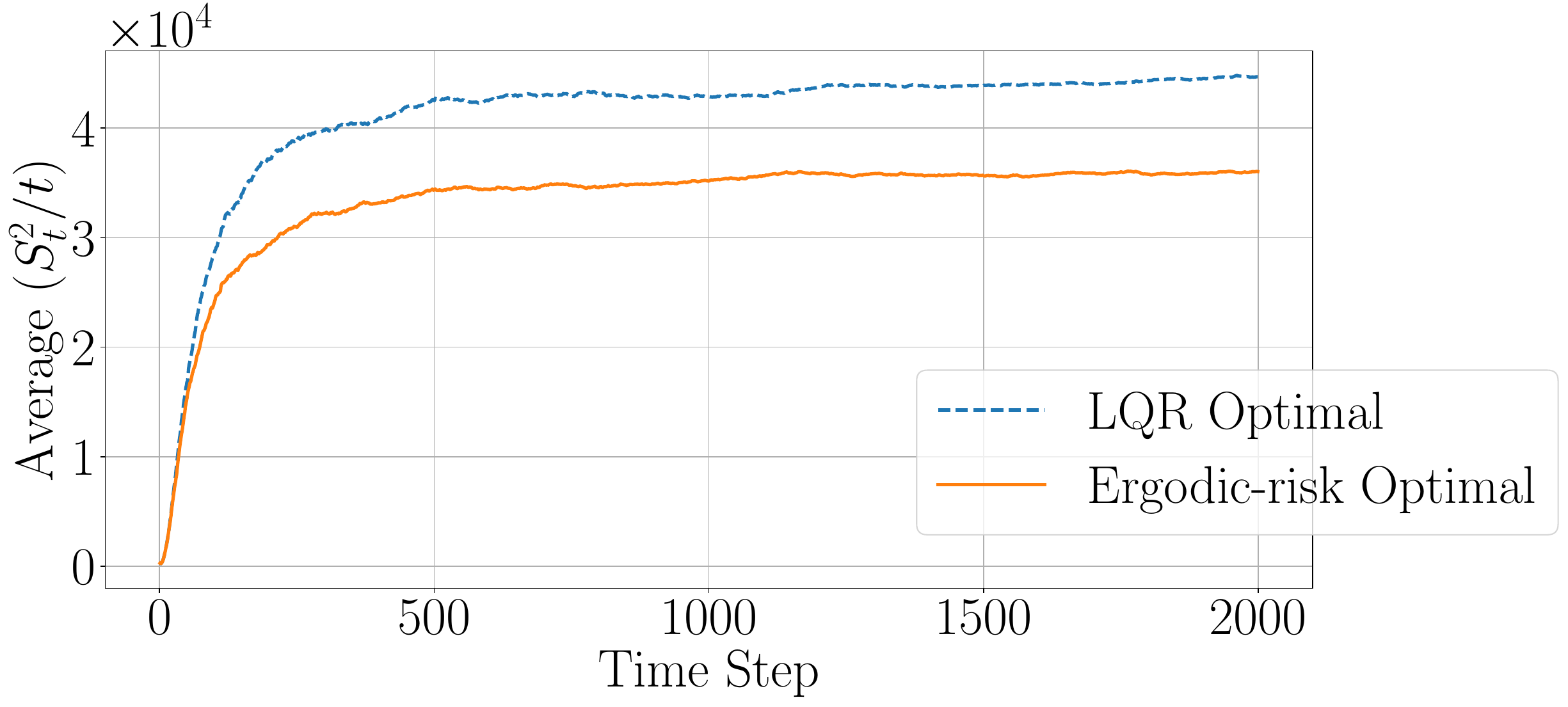}
        \caption{ The average of running covariance $S_t^2/t$  over 10,000 system responses under LQR versus the ergodic-risk optimal policies, approximating their corresponding asymptotic variance $\gamma_C^2$ as $t$ approaches infinity.}
    \label{fig:S-2-t}
    \vspace{-0.7cm}
\end{figure}

Finally, recall that $\gamma_N^2$ is an estimate of $\gamma_C^2 = \lim_{t\to \infty} \E[S_t^2/t]$, as discussed in \Cref{rem:asymp-var-conditional-var}. So, even though the optimal ergodic-risk policy directly bounds $\gamma_N^2$, we know that it should also result in lower asymptotic variance $\gamma_C^2$ compared to optimal LQR policy. In order to observe this numerically, we compare the performance of this policy in terms of the running variance of the cumulative risk, i.e. $\E [S_t^2/t]$, which is guaranteed to be well-defined if the controller is stabilizing and process noise has finite fourth moment \Cref{thm:tailored-clt}. 
In \Cref{fig:S-2-t}, we illustrate the average of $S_t^2/t$ with $S_t$ defined in \cref{eq:def-C-infty} over 10,000 roll-outs of the dynamics in \cref{eq:dynamics} with $\P_W$ being the Student's t-distribution with parameter $\nu =5$. As time goes to infinity, this numerically approximates the value of $\gamma_C^2$ for each policy. Thus, this simulation illustrates that the proposed ergodic-risk policy from \cref{eq:optimization-reform} results in smaller asymptotic (empirical) values of $\E [S_t^2/t]$, confirming a more risk averse behavior compared to the LQR optimal one in terms of $\gamma_C^2$ (which, in turn, coincides with the covariance of ergodic-risk criterion $C_\infty$ in \cref{eq:def-C-infty}).

\begin{figure}[pt]
  \centering    \includegraphics[width=0.8\textwidth]{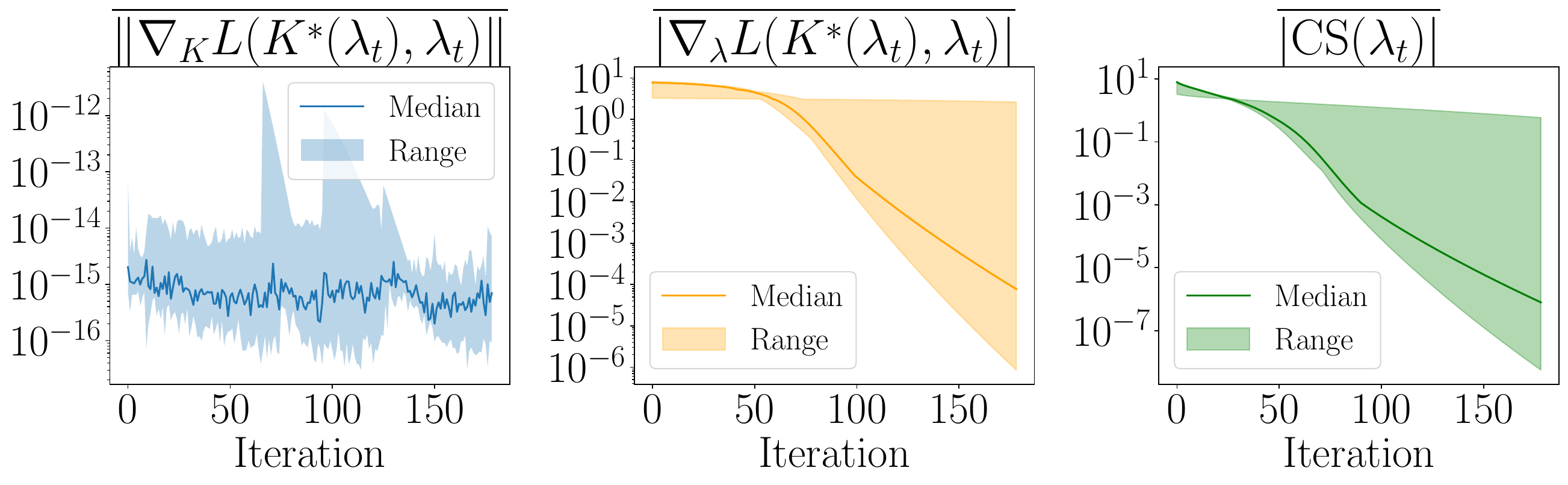}
  \caption{ Convergence of \Cref{algo} on 50 randomly sampled problem instances in terms of errors in KKT conditions vs iteration $m$.}
  \label{fig:mainfig}
  \vspace{-0.7cm}
\end{figure}

\textbf{Simulation 2:} Next, we illustrate the convergence behavior of \Cref{algo} on 50 randomly generated problem instances with 4 states and 2 inputs, where we set $\bar\beta$ such that any feasible policy is $10\%$ more conservative than the LQR optimal; i.e. $\bar\beta = 0.9 \gamma_N^2(K_{LQR})$. Its progress on 50 randomly sampled problem instances is illustrated in \Cref{fig:mainfig}. Because these instances are drawn at random, the Slater's condition might fail to hold, which essentially means such a conservative performance is not feasible by any policy $K$ and thus causing the algorithm to fail in such instances.%

\section{Conclusions and Future Directions}\label{sec:conclusion}
We introduced ergodic-risk criterion in COCP as a flexible framework to capture long-term cumulative uncertainties using appropriate risk functionals in presence of heavy-tailed disturbances. By considering linear constraints on $\gamma_N^2$ and integrating recent advances in policy optimization, we established strong duality and proposed a primal-dual algorithm with convergence guarantees. Key future directions include extending this approach to directly constrain $\gamma_C^2$, considering affine policy optimization, and developing efficient sample-based algorithms to obtain ergodic-risk policies directly from input-output data.

\section*{Acknowledgments}
This work was in part funded by NSF AI institute 2112085. The first author extends heartfelt gratitude to Professor Sean Meyn for the invaluable contribution in his seminal book \cite{meyn_markov_2009}, and the several insightful discussions regarding this work. The authors are also thankful to Dr. Derya Cansever for insightful discussions and feedback on the initial arXiv version of this work. Finally, the authors would like to thank the Associate Editor and Reviewers for their constructive comments and suggestions.

\bibliographystyle{siamplain}
\bibliography{SICON_LQR}

\begin{thebibliography}{10}

\bibitem{biswas_ergodic_2023}
{\sc A.~Biswas and V.~S. Borkar}, {\em Ergodic risk-sensitive control—{A}
  survey}, Annual Reviews in Control, 55 (2023), pp.~118--141,
  \url{https://doi.org/10.1016/j.arcontrol.2023.03.001},
  \url{https://www.sciencedirect.com/science/article/pii/S1367578823000068}
  (accessed 2024-09-16).

\bibitem{borkar_risk-constrained_2014}
{\sc V.~Borkar and R.~Jain}, {\em Risk-{Constrained} {Markov} {Decision}
  {Processes}}, IEEE Trans. Autom. Control,  (2014), pp.~2574--2579,
  \url{https://doi.org/10.1109/TAC.2014.2309262},
  \url{http://ieeexplore.ieee.org/document/6750726/} (accessed 2024-06-14).

\bibitem{borkar_controlled_1990}
{\sc V.~S. Borkar and M.~K. Ghosh}, {\em Controlled diffusions with
  constraints}, Journal of Mathematical Analysis and Applications, 152 (1990),
  pp.~88--108, \url{https://doi.org/10.1016/0022-247X(90)90094-V},
  \url{https://www.sciencedirect.com/science/article/pii/0022247X9090094V}
  (accessed 2025-07-07).

\bibitem{bosworth_linearized_1992}
{\sc J.~T. Bosworth}, {\em Linearized aerodynamic and control law models of the
  {X}-{29A} airplane and comparison with flight data}, vol.~4356 of {NASA}
  {Technical} {Memorandum}, National Aeronautics and Space Administration
  (NASA, Office of Management), 1992,
  \url{https://ntrs.nasa.gov/api/citations/19920009932/downloads/19920009932.pdf}.

\bibitem{chow_risk-constrained_2017}
{\sc Y.~Chow, M.~Ghavamzadeh, L.~Janson, and M.~Pavone}, {\em
  Risk-{Constrained} {Reinforcement} {Learning} with {Percentile} {Risk}
  {Criteria}}, Apr. 2017, \url{https://doi.org/10.48550/arXiv.1512.01629},
  \url{http://arxiv.org/abs/1512.01629} (accessed 2023-11-04).
\newblock arXiv:1512.01629 [cs, math].

\bibitem{durrett_probability_2019}
{\sc R.~Durrett}, {\em Probability: {Theory} and {Examples}}, Cambridge
  University Press, Cambridge, 5th~ed., 2019,
  \url{https://doi.org/10.1017/9781108591034},
  \url{https://www.cambridge.org/core/books/probability/DD9A1907F810BB14CCFF022CDFC5677A}
  (accessed 2024-03-12).

\bibitem{eichler_risks_2013}
{\sc H.-G. Eichler and B.~Bloechl-Daum, et.~al.}, {\em The risks of risk
  aversion in drug regulation}, Nature Reviews Drug Discovery, 12 (2013),
  pp.~907--916, \url{https://doi.org/10.1038/nrd4129},
  \url{https://www.nature.com/articles/nrd4129} (accessed 2024-09-13).

\bibitem{fazel_global_2018}
{\sc M.~Fazel, R.~Ge, S.~Kakade, and M.~Mesbahi}, {\em Global convergence of
  policy gradient methods for the linear quadratic regulator}, in Int. {Conf}.
  on {Machine} {Learning}, PMLR, July 2018, pp.~1467--1476,
  \url{http://proceedings.mlr.press/v80/fazel18a.html}.

\bibitem{hall_martingale_1980}
{\sc P.~Hall and C.~C. Heyde}, {\em Martingale {Limit} {Theory} and {Its}
  {Application}}, Academic Press, Dec. 1980,
  \url{https://shop.elsevier.com/books/martingale-limit-theory-and-its-application/birnbaum/978-0-12-319350-6}
  (accessed 2024-08-12).

\bibitem{kishida_risk-aware_2023}
{\sc M.~Kishida and A.~Cetinkaya}, {\em Risk-aware linear quadratic control
  using conditional value-at-risk}, IEEE Transactions on Automatic Control, 68
  (2023), pp.~416--423, \url{https://doi.org/10.1109/TAC.2022.3142131},
  \url{https://ieeexplore.ieee.org/document/9678012/} (accessed 2024-01-04).

\bibitem{komorowski_central_2012}
{\sc T.~Komorowski and A.~Walczuk}, {\em Central limit theorem for {Markov}
  processes with spectral gap in the {Wasserstein} metric}, Stochastic
  Processes and their Applications, 122 (2012), pp.~2155--2184,
  \url{https://doi.org/10.1016/j.spa.2012.03.006},
  \url{https://www.sciencedirect.com/science/article/pii/S0304414912000397}
  (accessed 2024-04-11).

\bibitem{kumar_stochastic_2015}
{\sc P.~R. Kumar and P.~P. Varaiya}, {\em Stochastic {Systems}: {Estimation},
  identification, and adaptive control}, Classics in applied mathematics,
  Society for Industrial and Applied Mathematics (SIAM, Philadelphia,PA),
  Philadelphia, Pennsylvania, 2015,
  \url{https://doi.org/10.1137/1.9781611974263} (accessed 2024-08-07).
\newblock OCLC: 930320873.

\bibitem{majumdar_how_2020}
{\sc A.~Majumdar and M.~Pavone}, {\em How {Should} a {Robot} {Assess} {Risk}?
  {Towards} an {Axiomatic} {Theory} of {Risk} in {Robotics}}, in Robotics
  {Research}, 2020, pp.~75--84,
  \url{https://doi.org/10.1007/978-3-030-28619-4_10}.

\bibitem{meyn_policy_1997}
{\sc S.~Meyn}, {\em The policy iteration algorithm for average reward {Markov}
  decision processes with general state space}, IEEE Transactions on Automatic
  Control, 42 (1997), pp.~1663--1680, \url{https://doi.org/10.1109/9.650016},
  \url{https://ieeexplore.ieee.org/document/650016} (accessed 2024-06-10).
\newblock Conference Name: IEEE Transactions on Automatic Control.

\bibitem{meyn_markov_2009}
{\sc S.~Meyn and R.~L. Tweedie}, {\em Markov {Chains} and {Stochastic}
  {Stability}}, Cambridge University Press, 2009,
  \url{https://link-springer-com.ezp-prod1.hul.harvard.edu/book/10.1007/978-1-4471-3267-7}
  (accessed 2024-04-24).

\bibitem{nedic_subgradient_2009}
{\sc A.~Nedić and A.~Ozdaglar}, {\em Subgradient methods for saddle-point
  problems}, J. Optim. Theory Appl., 142 (2009), pp.~205--228,
  \url{https://doi.org/10.1007/s10957-009-9522-7},
  \url{http://link.springer.com/10.1007/s10957-009-9522-7} (accessed
  2024-09-09).

\bibitem{rockafellar_optimization_2000}
{\sc R.~T. Rockafellar and S.~Uryasev}, {\em Optimization of conditional
  value-at-risk}, The Journal of Risk, 2 (2000), pp.~21--41,
  \url{https://doi.org/10.21314/JOR.2000.038},
  \url{http://www.risk.net/journal-of-risk/technical-paper/2161159/optimization-conditional-value-risk}
  (accessed 2023-11-04).

\bibitem{sopasakis_risk-averse_2019}
{\sc P.~Sopasakis, D.~Herceg, A.~Bemporad, and P.~Patrinos}, {\em Risk-averse
  model predictive control}, Automatica,  (2019), pp.~281--288,
  \url{https://doi.org/10.1016/j.automatica.2018.11.022},
  \url{https://www.sciencedirect.com/science/article/pii/S0005109818305545}
  (accessed 2024-09-12).

\bibitem{talebi_regularizability_2022}
{\sc S.~Talebi, S.~Alemzadeh, N.~Rahimi, and M.~Mesbahi}, {\em On
  {Regularizability} and its application to online control of unstable {LTI}
  systems}, IEEE Trans on Automatic Control, 67 (2022), pp.~6413--6428,
  \url{https://doi.org/10.1109/TAC.2021.3131148}.

\bibitem{talebi_policy_2023}
{\sc S.~Talebi and M.~Mesbahi}, {\em Policy optimization over submanifolds for
  linearly constrained feedback synthesis}, IEEE Transactions on Automatic
  Control,  (2023), pp.~1--16, \url{https://doi.org/10.1109/TAC.2023.3306384},
  \url{https://ieeexplore.ieee.org/document/10224332} (accessed 2024-03-26).

\bibitem{talebi_data-driven_2023}
{\sc S.~Talebi, A.~Taghvaei, and M.~Mesbahi}, {\em Data-driven optimal
  filtering for linear systems with unknown noise covariances}, in Advances in
  {Neural} {Inform}. {Process}. {Sys}., 2023, pp.~69546--69585,
  \url{https://proceedings.neurips.cc/paper_files/paper/2023/file/dbe8185809cb7032ec7ec6e365e3ed3b-Paper-Conference.pdf}.

\bibitem{talebi_policy_2024}
{\sc S.~Talebi, Y.~Zheng, S.~Kraisler, N.~Li, and M.~Mesbahi}, {\em Policy
  {Optimization} in {Control}: {Geometry} and {Algorithmic} {Implications}},
  June 2024, \url{https://doi.org/10.48550/arXiv.2406.04243},
  \url{http://arxiv.org/abs/2406.04243} (accessed 2025-07-30).
\newblock arXiv:2406.04243 [math].

\bibitem{tsiamis_risk-constrained_2020}
{\sc A.~Tsiamis, D.~S. Kalogerias, L.~F. Chamon, A.~Ribeiro, and G.~J. Pappas},
  {\em Risk-constrained linear-quadratic regulators}, in 59th {IEEE}
  {Conference} on {Decision} and {Control}, IEEE, 2020, pp.~3040--3047.

\bibitem{whittle_risk-sensitive_1981}
{\sc P.~Whittle}, {\em Risk-sensitive linear/quadratic/gaussian control},
  Advances in Applied Probability, 13 (1981), pp.~764--777,
  \url{https://doi.org/10.2307/1426972},
  \url{https://www.jstor.org/stable/1426972} (accessed 2024-03-01).

\bibitem{zhang_policy_2021}
{\sc K.~Zhang, B.~Hu, and T.~Başar}, {\em Policy {Optimization} for {H}-2
  {Linear} {Control} with {H}-infinity {Robustness} {Guarantee}: {Implicit}
  {Regularization} and {Global} {Convergence}}, SIAM Journal on Control and
  Optimization, 59 (2021), pp.~4081--4109,
  \url{https://doi.org/10.1137/20M1347942},
  \url{https://epubs.siam.org/doi/10.1137/20M1347942} (accessed 2024-03-07).

\bibitem{zhao_global_2023}
{\sc F.~Zhao, K.~You, and T.~Başar}, {\em Global convergence of policy
  gradient primal–dual methods for risk-constrained {LQRs}}, IEEE
  Transactions on Automatic Control, 68 (2023), pp.~2934--2949,
  \url{https://ieeexplore.ieee.org/abstract/document/10005813}.

\end{thebibliography}

\appendix

\section{Proofs of \Cref{lem:MDS} and \Cref{lem:ave-conv}}

\begin{proof}[Proof of \Cref{lem:MDS}]
    Note that $g(\x_t,\u_t) \in \F_t$ for each $t$. So, $C_t$ is $\F_{t}$-adapted.
    By Conditional Jensen's Inequality we obtain that for each $t$
    \begin{align*}
        \E|C_t| &\leq \E |g(\x_t,\u_t)| + \E[\E[|g(\x_t,\u_t)|\; |\F_{t-1}] ] \\
    &= 2 \E|g(\x_t,\u_t)| \\
    &\leq 2 M \left( \E\|\x_t\|^p+\E\|\u_t\|^p \right) <\infty,
    \end{align*}
    where the equality follows the Tower property and the last inequality is due to the moment condition on the process noise, the linear dynamics, and the fact that $\pi$ is affine. Because in that case, there exists another constant $M_2$ such that 
    \begin{align*}
        \|\x_{t}\| &\leq M_2(\|\x_{t-1}\| + \|\w_t\| + 1),\\
    \|\u_{t}\| &\leq M_2(\|\x_{t-1}\| + 1),
    \end{align*}
    which by H\"{o}lder inequality imply that
    \(\|\x_{t}\|^p + \|\u_{t}\|^p \leq 3^{p} M_2^p(\|\x_{t-1}\|^p + \|\w_t\|^p + 1).\)
    If the process noise and initial condition has moments up to order $p$ then $\E (\|\x_t\|^p + \|\u_t\|^p)$ is bounded for each finite $t$.
    Finally, $\E[C_t|\F_{t-1}] = 0$ by linearity of expectation and Tower property, and thus $\{C_t,\F_{t}\}$ is a MDS. 
\end{proof}

\begin{proof}[Proof of \Cref{lem:ave-conv}]
    Note that 
    \begin{equation}\label{eq:sys-traj}
        \x_{t+1} = A_K^{t+1} \x_0 + \bar{x}_t + \sum_{\tau =0}^{t} A_K^{t-\tau}H\w_{\tau+1}
    \end{equation}
    where $\bar{x}_t \coloneqq \sum_{\tau =0}^{t} A_K^\tau B \ell$.
    By \Cref{assmp:noise} on noise distribution, we observe that 
    \[\Sigma_{t+1} \coloneqq \E[(\x_{t+1}-\bar{x}_t) (\x_{t+1} - \bar{x}_t)^\intercal]\]
    must satisfy for all $t\geq0$:
    \begin{equation*}
        \Sigma_{t+1} =  A_K^{t+1}\Sigma_0 (A_K^\intercal)^{t+1} 
        + \sum_{\tau =0}^{t} A_K^{\tau} H \Sigma_W H^\intercal (A_K^\intercal)^{\tau}.
    \end{equation*}
    As $K\in\stableK$, $A_K$ is stabilizing and so $A_K^t$ converges to zero geometrically fast as $t$ goes to infinity. Thus, the series is absolutely convergent and we observe that $\Sigma_t \to \Sigma_K$ where the existence, uniqueness, and positive-definiteness of $\Sigma_K$ follows from Discrete Lyapunov Equation. Additionally, as $t\to \infty$ we obtain that $\bar{x}_t \to \bar{x}_K$, and thus $\E[\x_{t+1}] = A_K^{t+1}m_0 + \bar{x}_t \to \bar{x}_K$.
    Therefore, for all $t\geq1$, 
    \begin{equation}\label{eq:exp-Lambda-t}
        \E[\s_t] =  \sum_{s=1}^{t} A_K^s m_0 + \sum_{s=1}^{t} (\bar{x}_{s-1}-\bar{x}_K),
    \end{equation}
    and thus 
    $\E[\s_t/t] \to 0$.
    Note that by exponential stability, both
    \(\sum_{s=1}^t (\bar{x}_{s-1} - \bar{x}_K)\)
    and $\sum_{s=1}A_K^s m_0$ converge exponentially fast as $t\to \infty$. Thus, by \cref{eq:exp-Lambda-t}, $\E[\Lambda_t/\sqrt{t}] \to 0$.
    Similarly, 
    \begin{equation}\label{eq:expec-Gamma-t}
        \E[\Gamma_t] = \sum_{s=1}^t \Sigma_{s} 
        + \sum_{s=1}^t 2 A_K^s m_0 (\bar{x}_{s-1} - \bar{x}_K)^\intercal
        + (\bar{x}_{s-1} - \bar{x}_K) (\bar{x}_{s-1} - \bar{x}_K)^\intercal
    \end{equation}
    and thus
    $\E[\Gamma_t/t] \to \Sigma_K$ as $\Sigma_t \to \Sigma_K$.
    Next, as $\E[\x_t] \to \bar{x}_K$ exponentially fast, it suffices to show the last claim for the centered process $\{\x_t -\bar{x}_K\}$. So consider \cref{eq:sys-traj} where $\ell =0$ implying the $\bar{x}_K =0$ and thus
        \begin{align*}
            \s_{t+1} &= \sum_{s= 0}^t \x_{s+1}
            = \sum_{s=0}^t A_K^{s+1} \x_0 + \sum_{s=0}^t\sum_{\tau=0}^s A_K^{s-\tau} H \w_{\tau+1} \\
            &= - \x_0 + \sum_{s=-1}^{t} A_K^{s+1} \x_0 + \sum_{\tau=0}^t \sum_{s=\tau}^t A_K^{s-\tau} H \w_{\tau+1}= -\x_0 + \sum_{\tau=-1}^t \mathcal{A}_{t-\tau} H \w_{\tau+1},
        \end{align*}
        where $\mathcal{A}_{t} \coloneqq \sum_{s=0}^t A_K^s$ and for simplicity we denote, with abuse of notation, $H\w_{0} \coloneqq \x_0$. Let us without loss of generality assume $\Sigma_W = \Sigma_0$ in the rest of the argument.
        If $K\in\stableK$,  then there exists uniform constants $C,\rho$ such that \cite[Lemma 6]{talebi_data-driven_2023} $\|A_K^s\| \leq C \rho^s, \forall s$, implying that
        \[\|\mathcal{A}_t\| \leq C \frac{1-\rho^{t}}{1-\rho} \leq \frac{C}{1-\rho}\]
        Now, for any $\epsilon>0$, Markov Inequality implies
        \begin{align*}
            \P\{\|\s_{t+1}+\x_0\| > \epsilon t\}
            =& \P\{\tr{(\s_{t+1}+\x_0)(\s_{t+1}+\x_0)^\intercal} > \epsilon^2 t^2\}\\
            \leq & \frac{1}{\epsilon^2 t^2} \E \tr{(\s_{t+1}+\x_0)(\s_{t+1}+\x_0)^\intercal} \\
            = &  \frac{1}{\epsilon^2 t^2} \E \tr{\sum_{\tau=-1}^t \sum_{s=-1}^t \mathcal{A}_{t-\tau} H \w_\tau \w_s^\intercal H^\intercal \mathcal{A}_{t-s}^\intercal}\\
            = &  \frac{1}{\epsilon^2 t^2} \tr{ \sum_{s=-1}^t \mathcal{A}_{t-s} H \Sigma_W H^\intercal \mathcal{A}_{t-s}^\intercal}\\
            = &  \frac{1}{\epsilon^2 t^2} \tr{ H\Sigma_W H^\intercal \sum_{s=0}^{t+1} \mathcal{A}_{s}^\intercal \mathcal{A}_{s}}\\
            \leq &  \frac{\|H\|^2\tr{ \Sigma_W}}{\epsilon^2 t^2}  \left\|\sum_{s=0}^{t+1} \mathcal{A}_{s}^\intercal \mathcal{A}_{s}\right\| \\
            \leq &  \frac{\|H\|^2\tr{ \Sigma_W}}{\epsilon^2 t^2}  \sum_{s=0}^{t+1} \left\|\mathcal{A}_{s}\right\|^2\\
            \leq & \frac{\|H\|^2\tr{ \Sigma_W}}{\epsilon^2 t^2} \frac{C^2 (t+1)}{(1-\rho)^2} 
            \leq \mathcal{O}(\frac{1}{\epsilon^2 t}),
        \end{align*}
        where we used $\E{\w_\tau\w_s^\intercal} =0$ for $\tau\neq s$ and bounded second-order moments of the noise and initial condition. This implies the convergence of $\s_t/t \xrightarrow{p} 0 $ in probability.
        For the last claim, the case of $\ell = 0$ follows by a similar argument using Markov Inequality and convergence of $\Sigma_t$ and thus, is omitted. For the case of nonzero $\ell$, by \cite[Lemma 6]{talebi_policy_2023}, there exists a constant $\bar c_2$ such that $\|\bar{x}_t\|\leq \bar c_2 \|\ell\|$ for all $t\geq 0$. Therefore, $\|\x_t-\bar{x}_t\|\geq \|\x_t\| - \bar c_2\|\ell\|$ which implies that
        \[\P \left\{ \|\x_t\|\geq m \right\} \leq \P\left\{ \|\x_t - \bar{x}_t\|\geq m - \bar c_2\|\ell\|\right\}.\]
    This completes the proof. 
\end{proof}

\end{document}